# HOW MANY ENTRIES OF A TYPICAL ORTHOGONAL MATRIX CAN BE APPROXIMATED BY INDEPENDENT NORMALS?[1]

By Tiefeng Jiang

*University of Minnesota*

We solve an open problem of Diaconis that asks what are the largest orders of $p_n$ and $q_n$ such that $Z_n$, the $p_n \times q_n$ upper left block of a random matrix $\mathbf{\Gamma}_n$ which is uniformly distributed on the orthogonal group $O(n)$, can be approximated by independent standard normals? This problem is solved by two different approximation methods.

First, we show that the *variation distance* between the joint distribution of entries of $Z_n$ and that of $p_n q_n$ independent standard normals goes to zero provided $p_n = o(\sqrt{n})$ and $q_n = o(\sqrt{n})$. We also show that the above variation distance does not go to zero if $p_n = [x\sqrt{n}]$ and $q_n = [y\sqrt{n}]$ for any positive numbers $x$ and $y$. This says that the largest orders of $p_n$ and $q_n$ are $o(n^{1/2})$ in the sense of the above approximation.

Second, suppose $\mathbf{\Gamma}_n = (\gamma_{ij})_{n \times n}$ is generated by performing the Gram–Schmidt algorithm on the columns of $\mathbf{Y}_n = (y_{ij})_{n \times n}$, where $\{y_{ij}; 1 \leq i, j \leq n\}$ are i.i.d. standard normals. We show that $\varepsilon_n(m) := \max_{1 \leq i \leq n, 1 \leq j \leq m} |\sqrt{n} \cdot \gamma_{ij} - y_{ij}|$ *goes to zero in probability* as long as $m = m_n = o(n/\log n)$. We also prove that $\varepsilon_n(m_n) \to 2\sqrt{\alpha}$ in probability when $m_n = [n\alpha/\log n]$ for any $\alpha > 0$. This says that $m_n = o(n/\log n)$ is the largest order such that the entries of the first $m_n$ columns of $\mathbf{\Gamma}_n$ can be approximated simultaneously by independent standard normals.

**1. Introduction.** Let $\mathbf{\Gamma}_n = (\gamma_{ij})$ be a random orthogonal matrix which is uniformly distributed on the orthogonal group $O(n)$. Let $Z_n$ be the $p_n \times q_n$ upper left block of $\mathbf{\Gamma}_n$, where $p_n$ and $q_n$ are two positive integers. The open problem in Section 6.3 from [10] is as follows: what are the largest orders of $p_n$ and $q_n$ such that the variation distance between the joint distribution

Received February 2003; revised August 2005.
[1]Supported in part by NSF Grants DMS-03-08151 and DMS-04-49365.

of the entries of $Z_n$ and that of $p_n q_n$ independent standard normals goes to zero as $n \to \infty$. We answer this question here. Before stating the results formally, let us first review some history of this problem.

In studying "Equivalence of Ensembles" in statistical mechanics, Borel [5] showed that

$$(1.1) \qquad P(\sqrt{n}\gamma_{11} \leq x) \to \frac{1}{\sqrt{2\pi}} \int_{-\infty}^{x} e^{-t^2/2}\,dt$$

as $n \to \infty$ for any real number $x$. For more information about this formula, one is referred to [27], page 197 in [26], page 412 in [23], page 342 in [4], [11] and [24, 25].

Similar results for fixed $m$ are derived through Brownian motion by Gallardo [16] and Yor [32]. Let $\boldsymbol{\gamma}_1$ be the first column of $\boldsymbol{\Gamma}_n$. Stam [30] proved that $d_m$, the variation distance between the distribution of the first $m$ coordinates of $\boldsymbol{\gamma}_1$ and the distribution of $m$ independent standard normals, goes to zero provided $m = o(\sqrt{n})$ as $n \to \infty$. He applied this result to a geometric probability problem.

In studying a finite representation theorem of the de Finetti type, Diaconis and Freedman [11] showed that the above $d_m$ goes to zero as $n \to \infty$ provided $m = o(n)$. On the other hand, in studying a de Finetti-type theorem on a finite sequence of orthogonal invariant random vectors, Diaconis, Eaton and Lauritzen [14] proved the following.

THEOREM A.1. *For each $n \geq 1$, let $Z_n$ be the $p_n \times q_n$ upper left block of a random matrix $\boldsymbol{\Gamma}_n$ which is uniformly distributed on the orthogonal group $O(n)$. Let also $\delta_n$ be the variation distance between the distribution of the $p_n q_n$ entries of $Z_n$ and the joint distribution of $p_n q_n$ independent standard normals. Then $\delta_n \to 0$ if $p_n = o(n^\alpha)$ and $q_n = o(n^\alpha)$ for $\alpha = 1/3$.*

Since the publication of [14], there have been various speculations on the maximum value $\alpha$ to make the variation distance go to zero. Here are three major ones: (a) $p_n = O(n^{1/3})$ and $q_n = O(n^{1/3})$; (b) $p_n = o(n^{1/2})$ and $q_n = o(n^{1/2})$; (c) $p_n = o(n)$ and $q_n = o(n)$. Recently Collins [7] showed that the variation distance in Theorem A.1 goes to zero when $p_n = O(n^{1/3})$ and $q_n = O(n^{1/3})$.

Attempts to improve on the orders of $p_n$ and $q_n$ are partly motivated by the following reasons. First, it is well known that the above $\boldsymbol{\Gamma}_n$ is close to $\boldsymbol{\Gamma}'_n$, an $n \times n$ matrix with independent normals as entries. Mathematically, it is interesting to know in what sense they are close. Diaconis and Shahshahani [12], Diaconis and Evans [13] and Rains [28] characterized relationships between the traces of $\boldsymbol{\Gamma}_n$ and those of $\boldsymbol{\Gamma}'_n$ in terms of expectations; Johansson [20] obtained the speed of convergence of traces of $\boldsymbol{\Gamma}_n$ to a normal random variable; D'Aristotile, Diaconis and Newman [8] showed that the linear combination of entries of $\mathbf{\Gamma}_n$ also converges weakly to a normal distribution. Second, improving the orders of $p_n$ and $q_n$ has a lot of applications; see [14] and [19]. In the last paper Jiang also proved the following coupling result.

THEOREM A.2. *For each $n \geq 2$, there exists matrices $\mathbf{\Gamma}_n = (\gamma_{ij})_{1 \leq i,j \leq n}$ and $\mathbf{\Gamma}'_n = (\gamma'_{ij})_{1 \leq i,j \leq n}$ whose $2n^2$ elements are random variables defined on the same probability space such that:*

(i) *the law of $\mathbf{\Gamma}_n$ is the normalized Haar measure on the orthogonal group $O_n$;*

(ii) *$\{\gamma'_{ij}; 1 \leq i, j \leq n\}$ are independent standard normals;*

(iii) *set $\varepsilon_n(m) = \max_{1 \leq i \leq n, 1 \leq j \leq m} |\sqrt{n}\gamma_{ij} - \gamma'_{ij}|$ for $m = 1, 2, \ldots, n$. Then*

$$\varepsilon_n(m_n) \to 0 \quad \text{in probability,}$$

*as $n \to \infty$ provided $m_n = o(n/(\log n)^2)$.*

It says that $n^2/(\log n)^2$ elements of $\mathbf{\Gamma}_n$ can be approximated by the corresponding elements of $\mathbf{\Gamma}'_n$ in terms of convergence in probability, which is weaker than the convergence in variation norm.

This theorem highlights the interest in improving the orders of $p_n$ and $q_n$. It seems to suggest that Theorem A.1 holds for much larger $p_n$ and $q_n$. This is why people conjectured that the maximum orders of $p_n$ and $q_n$ are $o(n)$. At the same time it would be interesting to know the largest order of $m_n$ such that Theorem A.2 holds.

In this paper we prove that the maximum value of $\alpha$ as in Theorem A.1 is actually $1/2$, and the largest order of $m_n$ such that $\varepsilon_n(m_n) \to 0$ in probability is $o(n/\log n)$, where $\varepsilon_n(m_n)$ is as in Theorem A.2. To state our results formally, let us recall the definition of variation distance first.

Let $\mu$ and $\nu$ be two probability measures on $(\mathbb{R}^m, \mathcal{B})$, where $\mathcal{B}$ is the Borel $\sigma$-algebra. The variation distance between $\mu$ and $\nu$, denoted by $\|\mu - \nu\|$, is equal to

$$(1.2) \quad \|\mu - \nu\| = 2 \cdot \sup_{A \in \mathcal{B}} |\mu(A) - \nu(A)| = \int_{\mathbb{R}^m} |f(x) - g(x)| \, dx_1 \, dx_2 \cdots dx_m,$$

provided $\mu$ and $\nu$ have density functions $f(x)$ and $g(x)$ with respect to the Lebsegue measure, respectively. For each $n \geq 1$, suppose that $Z_n$ is the $p_n \times q_n$ upper left block of a random matrix $\mathbf{\Gamma}_n$ which is uniformly distributed on the orthogonal group $O(n)$. Let $G_n$ be the joint distribution of $p_n q_n$ independent standard normals. We use $\mathcal{L}(\sqrt{n} Z_n)$ to represent the joint probability distribution of the $p_n q_n$ random entries of $\sqrt{n} Z_n$. It is not difficult to see that $\|\mathcal{L}(\sqrt{n} Z_n) - G_n\|$ is nondecreasing in $p_n$ and $q_n$, respectively.

THEOREM 1. *If $p_n = o(\sqrt{n})$ and $q_n = o(\sqrt{n})$ as $n \to \infty$, then*

$$\lim_{n \to \infty} \|\mathcal{L}(\sqrt{n} Z_n) - G_n\| = 0.$$

As usual, the notation $[a]$ stands for the integer part of a positive integer $a$.

THEOREM 2. *Let $x > 0$ and $y > 0$ be two numbers and $p_n = [xn^{1/2}]$ and $q_n = [yn^{1/2}]$. Then*

$$\liminf_{n \to \infty} \|\mathcal{L}(\sqrt{n}Z_n) - G_n\| \geq \phi(x, y) > 0,$$

*where $\phi(x,y) := E|\exp(-\frac{x^2 y^2}{8} + \frac{xy}{4}\xi) - 1| \in (0, 1)$ and $\xi$ is a standard normal.*

One can see that $\phi(0,0) = 0$, which roughly reflects the flavor of Theorem 1. This is rigorously true if the conclusion in Theorem 2 is replaced by $\lim_{n \to \infty} \|\mathcal{L}(\sqrt{n}Z_n) - G_n\| = \phi(x, y)$. A further analysis shows that the inequality in the theorem is actually strict.

Why are the maximum orders of $p_n$ and $q_n$ equal to $o(n^{1/2})$ as shown in Theorems 1 and 2?

There are two reasons. First, Diaconis and Freedman [11] showed that the variation distance between the distribution of the $o(n)$ entries of the first column of $\mathbf{\Gamma}_n$ and that of independent normals goes to zero. We know that $Z_n$, a $p_n$ by $q_n$ sub-matrix of $\mathbf{\Gamma}_n$, has $p_n q_n$ elements. One can guess that the number of approximated entries are fixed (loosely speaking). So the largest $\alpha$ in $p_n = o(n^\alpha)$ and $q_n = o(n^\alpha)$ must be $1/2$. Second, we can see this mathematically. Let $f_n(z)$ and $g_n(z)$ be the density functions of $\sqrt{n}Z_n$ and $G_n$, respectively. By (1.2),

$$(1.3) \quad \|\mathcal{L}(\sqrt{n}Z_n) - G_n\| = \int \left|\frac{f_n(z)}{g_n(z)} - 1\right| g_n(z)\,dz = E\left|\frac{f_n(X_n)}{g_n(X_n)} - 1\right|,$$

where the integration region in the first integral is $\mathbb{R}^{p_n q_n}$, and the $p_n q_n$ entries of the matrix $X_n$ are independent standard normals. The term $f(X_n)/g(X_n)$, as will be shown later, converges weakly to a lognormal distribution when both $p_n$ and $q_n$ are of order $n^{1/2}$; $f(X_n)/g(X_n)$ converges to one when both $p_n$ and $q_n$ are of order $o(n^{1/2})$.

Now we consider the approximation method as in Theorem A.2.

Let $\mathbf{Y}_n = (y_{ij})_{1 \leq i,j \leq n}$, where $y_{ij}$'s are independent standard normals. Let also $\mathbf{\Gamma}_n = (\gamma_{ij})_{1 \leq i,j \leq n}$ be the orthogonal matrix obtained from performing the Gram–Schmidt procedure on the columns of $\mathbf{Y}_n$ (the procedure is briefly reviewed at the beginning of Section 3). Define

$$\varepsilon_n(m) = \max_{1 \leq i \leq n, 1 \leq j \leq m} |\sqrt{n}\gamma_{ij} - y_{ij}|.$$

We have the following theorem.

THEOREM 3. *Let $\{m_n < n;\ n \geq 1\}$ be a sequence of positive integers. Then:*

(i) the matrix $\mathbf{\Gamma}_n$ is Haar invariant on the orthogonal group $O(n)$;
(ii) $\varepsilon_n(m_n) \to 0$ in probability, provided $m_n = o(n/\log n)$ as $n \to \infty$;
(iii) for any $\alpha > 0$, we have that $\varepsilon_n([n\alpha/\log n]) \to 2\sqrt{\alpha}$ in probability as $n \to \infty$.

This theorem tells us that the maximum order of $m_n$ such that $\varepsilon_n(m_n) \to 0$ in probability is that $m_n = o(n/\log n)$, where the typical orthogonal matrix $\mathbf{\Gamma}_n$ is obtained through performing the Gram–Schmidt procedure for a matrix whose elements are independent standard normals.

We prove Theorems 1 and 2 in Section 2. Theorem 3 is proved in Section 3. Technical lemmas used in Sections 2 and 3 are given in Section 4. At last, a couple of known results needed for the proof of Theorem 3 are listed in the Appendix.

**2. The proofs of Theorems 1 and 2.** First we list some lemmas needed for the proofs of Theorems 1 and 2. The proofs of these lemmas are listed in Section 4.1.

LEMMA 2.1. *Let $\Gamma(x)$, $x > 0$ be the standard Gamma function. Then:*

(i) $\quad 1 - \dfrac{1}{6n} < \dfrac{\Gamma(n + (1/2))}{\sqrt{n}\,\Gamma(n)} < 1 \qquad$ *for all $n \geq 1$;*

(ii) $\quad \left| \dfrac{\Gamma((n+1)/2)}{\sqrt{n/2}\,\Gamma(n/2)} - 1 \right| < \dfrac{3}{5n} \qquad$ *for all $n \geq 1$.*

LEMMA 2.2. *Let $f(u,v)$ be a real-valued function. Suppose the three second-order derivatives of $f$ exist, bounded below and above by $-M$ and $M$, respectively, over $[a,b] \times [c,d]$. Then*

$$\frac{1}{n^2} \sum_{j=j_1}^{j_2} \sum_{i=i_1}^{i_2} f\left(\frac{i}{n}, \frac{j}{n}\right) = \int_{j_1/n}^{(j_2+1)/n} \int_{i_1/n}^{(i_2+1)/n} f(x,y)\,dx\,dy$$
$$- \frac{1}{2n^3} \sum_{j=j_1}^{j_2} \sum_{i=i_1}^{i_2} f'_x\left(\frac{i}{n}, \frac{j}{n}\right)$$
$$- \frac{1}{2n^3} \sum_{j=j_1}^{j_2} \sum_{i=i_1}^{i_2} f'_y\left(\frac{i}{n}, \frac{j}{n}\right) + \varepsilon,$$

*where $|\varepsilon| \leq (i_2 - i_1)(j_2 - j_1)M/n^4$ for any $i_1, i_2, j_1$ and $j_2$ such that $na \leq i_1 < i_2 \leq nb - 1$ and $nc \leq j_1 < j_2 \leq nd - 1$.*

We will use the following setting a couple of times.

(2.1)   Let $X = (x_{ij})$ be a $p$ by $q$ matrix, where $\{x_{ij},\ 1 \leq i \leq p; 1 \leq j \leq q\}$ are i.i.d. standard normals. Let $\lambda_1, \lambda_2, \ldots, \lambda_q$ be the eigenvalues of $X'X$.

A sequence $\{X_n; n \geq 1\}$ will be studied, where $X_n$ is of the above setting for each $n$. We still use notation $X$ for $X_n$ sometimes when there is no confusion.

The next lemma is a standard result when using the moment method to show weak convergence of certain functions of eigenvalues of matrices with independent and identically distributed random variables as entries. It is can be seen from, for example, (2.15) and (2.16) in [3].

LEMMA 2.3.   *Let $\{p_n; n \geq 1\}$ and $\{q_n; n \geq 1\}$ be two sequences of positive integers such that $p_n \to \infty$ and $p_n/q_n \to \eta \in (0, \infty)$. For each $n$, assume the setting in (2.1) with $p = p_n$ and $q = q_n$. The following two statements hold. For each integer $k \geq 1$,*

$$\text{(i)} \quad E(\operatorname{tr}(X_n' X_n)^k) \sim p_n^k q_n \sum_{r=0}^{k-1} \frac{1}{r+1} \left(\frac{q_n}{p_n}\right)^r \binom{k}{r} \binom{k-1}{r}$$

*as $n \to \infty$.*

$$\text{(ii)} \quad \frac{\operatorname{tr}((X_n' X_n)^k)}{q_n^{k+1}} \to \sum_{r=0}^{k-1} \frac{\eta^{k-r}}{r+1} \binom{k}{r} \binom{k-1}{r}$$

*in probability as $n \to \infty$.*

LEMMA 2.4.   *Let $\varepsilon \in (0, 1)$. Let $\{p_n; n \geq 1\}$ and $\{q_n; n \geq 1\}$ be two sequences of positive integers such that $\varepsilon \leq p_n/q_n \leq \varepsilon^{-1}$ for all $n \geq 1$. For each $n$, assume the setting in (2.1) with $p = p_n$ and $q = q_n$. Assume $p_n \to \infty$ as $n \to \infty$. Then:*

(i)  $\operatorname{Var}(\operatorname{tr}((X_n' X_n)^2)) \sim p_n^2 q_n^2 + 8 p_n q_n (p_n + q_n)^2$ *as $n \to \infty$;*
(ii) $\operatorname{Cov}(\operatorname{tr}(X_n' X_n), \operatorname{tr}((X_n' X_n)^2)) \sim 4 p_n q_n (p_n + q_n)$ *as $n \to \infty$.*

The following lemma is Proposition 2.1 from [14] or Proposition 7.3 from [15]. This is the starting point of the proofs of Theorems 1 and 2.

LEMMA 2.5.   *Let $U$ be an $n$ by $n$ random matrix which is uniformly distributed on the orthogonal group $O_n$ and let $Z$ be the upper left $p \times q$ corner block of $U$. If $p + q \leq n$ and $q \leq p$, then the joint density function of entries of $Z$ is*

$$(2.2) \quad f(z) = (\sqrt{2\pi})^{-pq} \frac{\omega(n-p, q)}{\omega(n, q)} \{\det(I_q - z'z)^{(n-p-q-1)/2}\} I_0(z'z),$$

where $I_0(z'z)$ is the indicator function of the set that all $q$ eigenvalues of $z'z$ are in $(0,1)$, and $\omega(\cdot,\cdot)$ is the Wishart constant defined by

$$\frac{1}{\omega(r,s)} = \pi^{s(s-1)/4} 2^{rs/2} \prod_{j=1}^{s} \Gamma\left(\frac{r-j+1}{2}\right).$$

Here $s$ is a positive integer and $r$ is a real number, $r > s - 1$. When $p \leq q$, the density of $Z$ is obtained by interchanging $p$ and $q$ in the above Wishart constant.

To simplify notation, when there is no confusion, we write $p$ for $p_n$ and $q$ for $q_n$.

Let $g(z)$ be the joint density function of entries of $X = (x_{ij})_{p \times q}$, where $x_{ij}$'s are independent standard normals. So, $g(z) = (2\pi)^{-pq/2} \exp(-\operatorname{tr}(z'z)/2)$, where $z$ is a $p$ by $q$ matrix. We need to understand the ratio $f(z)/g(z)$ in later proofs. Assuming the $pq$ entries of $z$ are independent standard normals, then $f(z)/g(z)$ can be written as a product of a constant part and a random part. They are analyzed in the following two lemmas.

LEMMA 2.6. *Given $x > 0$ and $y > 0$, let $p = p_n = [xn^{1/2}]$ and $q = q_n = [yn^{1/2}]$. Set*

$$K_n = \left(\frac{2}{n}\right)^{pq/2} \prod_{j=1}^{q} \frac{\Gamma((n-j+1)/2)}{\Gamma((n-p-j+1)/2)}.$$

*Then*

(2.3) $$K_n = \exp\left\{-\left(\frac{p^2 q + pq^2}{4n} + \frac{xy}{4} + \frac{2x^3 y + 2xy^3 + 3x^2 y^2}{24}\right) + o(1)\right\},$$

*as $n$ is sufficiently large.*

PROOF. Suppose $p = 2k$. Using the fact that $\Gamma(x+1) = x\Gamma(x)$, we have that

(2.4)
$$K_n = \left(\frac{2}{n}\right)^{pq/2} \prod_{j=1}^{q}\left\{\left(\frac{2}{n}\right)^{-p/2} \prod_{i=1}^{k}\left(1 - \frac{2i+j-1}{n}\right)\right\}$$
$$= \prod_{j=1}^{q}\prod_{i=1}^{k}\left(1 - \frac{2i+j-1}{n}\right) := e^{B_n},$$

where

$$B_n := \sum_{j=0}^{q-1}\sum_{i=1}^{k} \log\left(1 - \frac{2i+j}{n}\right).$$

Let $f(s,t) = \log(1 - 2s - t)$ with $2s + t < 1$. Then $f'_s(s,t) = -2/(1-2s-t) = -2 + O(n^{-1/2})$, $f'_t(s,t) = -1/(1-2s-t) = -1 + O(n^{-1/2})$ and

$$\left|\frac{\partial^2 f}{\partial s^2}\right| = \left|\frac{-4}{(1-2s-t)^2}\right| \leq 5, \qquad \left|\frac{\partial^2 f}{\partial t^2}\right| = \left|\frac{-1}{(1-2s-t)^2}\right| \leq 5$$

and

$$\left|\frac{\partial^2 f}{\partial s \, \partial t}\right| = \left|\frac{-2}{(1-2s-t)^2}\right| \leq 5$$

for all $(s,t) \in [0, p/n] \times [0, q/n]$, as $n$ is sufficiently large. By Lemma 2.2,

$$B_n = n^2 \int_0^{q/n} \int_{1/n}^{(k+1)/n} \log(1 - 2s - t) \, ds \, dt + \frac{3kq}{2n} + O\left(\frac{1}{\sqrt{n}}\right)$$

(2.5)
$$= \frac{n^2}{2} \int_0^v \int_0^u \log(1 + s + t) \, ds \, dt$$
$$- \frac{n^2}{2} \int_0^v \int_0^{-2/n} \log(1 + s + t) \, ds \, dt + \frac{3xy}{4} + O\left(\frac{1}{\sqrt{n}}\right),$$

as $n$ is sufficiently large, where $u = -(p+2)/n$, $v = -q/n$. We now estimate the above integral. By Taylor's expansion, there exists $\delta > 0$ such that

$$\left|\log(1 + s + t) - \left((s+t) - \frac{(s+t)^2}{2}\right)\right| \leq (s+t)^3$$

for all $s$ and $t$ such that $s + t \in (0, \delta)$. Thus,

(2.6)
$$\left|\int_0^v \int_0^u \log(1 + s + t) \, ds \, dt - \left[\int_0^v \int_0^u (s+t) \, ds \, dt - \frac{1}{2} \int_0^v \int_0^u (s+t)^2 \, ds \, dt\right]\right|$$
$$\leq \int_0^v \int_0^u (s+t)^3 \, ds \, dt,$$

as both $u$ and $v$ are in $(0, \delta/2)$. It is trivial to verify that

$$\int_0^v \int_0^u (s+t)^k \, ds \, dt = \frac{1}{(k+1)(k+2)}((u+v)^{k+2} - u^{k+2} - v^{k+2})$$

for $k \geq 0$. Plugging this into (2.6), we obtain

$$\int_0^v \int_0^u \log(1 + s + t) \, ds \, dt = \left[\frac{u^2 v + uv^2}{2} - \frac{1}{12}(2uv^3 + 2u^3 v + 3u^2 v^2)\right]$$
$$+ O((u+v)^5),$$

as $n \to \infty$. [The actual formula for the integral is

$$\int_0^v \int_0^u \log(1+s+t)\,ds\,dt$$
$$= \tfrac{1}{2}(1+u+v)^2 \log(1+u+v)$$
$$\quad - \tfrac{1}{2}(1+u)^2 \log(1+u) - \tfrac{1}{2}(1+v)^2 \log(1+v) - \tfrac{3}{2}uv.]$$

Now substituting $u = -(p+2)/n$ and $v = -q/n$ back into the two integrals in (2.5), we have that

$$\frac{n^2}{2} \int_0^v \int_0^u \log(1+s+t)\,ds\,dt$$
(2.7)
$$= -\left[\frac{p^2 q + pq^2}{4n} + xy + \frac{y^2}{2} + \frac{2xy^3 + 2x^3 y + 3x^2 y^2}{24}\right] + O\!\left(\frac{1}{\sqrt{n}}\right)$$

and

(2.8) $\qquad \dfrac{n^2}{2} \displaystyle\int_0^v \int_0^{-2/n} \log(1+s+t)\,ds\,dt = -\dfrac{y^2}{2} + O\!\left(\dfrac{1}{\sqrt{n}}\right),$

as $n$ is sufficiently large. Combining (2.4), (2.5), (2.7) and (2.8), we obtain

(2.9) $\quad K_n = \exp\!\left\{-\left(\dfrac{p^2 q + pq^2}{4n} + \dfrac{xy}{4} + \dfrac{2x^3 y + 2xy^3 + 3x^2 y^2}{24}\right) + O(n^{-1/2})\right\},$

as $n$ is sufficiently large.

Now, suppose $p = 2k - 1$. Let

$$C_n = \prod_{j=1}^q \frac{\Gamma((n-j-p+1)/2)}{\Gamma((n-j-p)/2)\sqrt{(n-j-p)/2}}.$$

By Lemma 2.1, the $j$th term in the product, say, $C_{n,j}$, has the following property:

$$1 - \frac{1}{n-p-q} \leq C_{n,j} \leq 1 + \frac{1}{n-p-q}$$

for all $j = 1, 2, \ldots, q$ as long as $p + q \leq n - 3$. Therefore,

$$\left(1 - \frac{1}{n-p-q}\right)^q \leq C_n \leq \left(1 + \frac{1}{n-p-q}\right)^q.$$

Since $(1+x_n)^{k_n} = 1 + O(k_n x_n)$ as $x_n \to 0$, $k_n \to \infty$ and $k_n x_n \to 0$. It follows that $C_n = 1 + O(n^{-1/2})$, provided $p = O(\sqrt{n})$ and $q = O(\sqrt{n})$. So

$$K_n = \frac{1}{C_n}\left(\frac{2}{n}\right)^{pq} \prod_{j=1}^q \frac{\Gamma((n-j+1)/2)}{\Gamma((n-j-2k+1)/2)\sqrt{(n-j-2k+1)/2}}$$

$$\sim \left\{\prod_{j=1}^q \prod_{i=1}^k \frac{n-2i-j+1}{n}\right\} \cdot \left\{\prod_{j=1}^q \frac{n-j-2k+1}{n}\right\}^{-1/2} := K'_n \cdot K''_n,$$

where the fact $\Gamma(x+1) = x\Gamma(x)$ is used in the second step. Now

$$\log K_n'' = -\frac{1}{2}\sum_{j=1}^{q}\log\left(1 - \frac{j+2k-1}{n}\right)$$

(2.10)
$$= \frac{1}{2n}\sum_{j=1}^{q}(j+2k-1) + O\left(\frac{1}{\sqrt{n}}\right)$$

$$= \frac{y^2 + 2xy}{4} + O\left(\frac{1}{\sqrt{n}}\right),$$

as $n \to \infty$. In notation, $K_n'$ is identical to $K_n$ in (2.4). Keep in mind that the $k$ in (2.4) is equal to $p/2$; but the $k$ in the definition of $K_n'$ is equal to $(p+1)/2$. Apply (2.9) to $K_n'$ to obtain

$$-\log K_n' = \frac{(p+1)^2 q + (p+1)q^2}{4n} + \frac{xy}{4} + \frac{2x^3 y + 2xy^3 + 3x^2 y^2}{24} + O(n^{-1/2})$$

$$= \frac{p^2 q + pq^2}{4n} + \frac{3xy + y^2}{4} + \frac{2x^3 y + 2xy^3 + 3x^2 y^2}{24} + O(n^{-1/2}).$$

This together with (2.10) thus yields (2.3). □

LEMMA 2.7. *Suppose $x > 0$ and $y > 0$. For each $n \geq 1$, assume the setting in (2.1) with $p = p_n = [x\sqrt{n}]$ and $q = q_n = [y\sqrt{n}]$. Define*

$$L_n = \left\{\prod_{i=1}^{q}\left(1 - \frac{\lambda_i}{n}\right)\right\}^{(n-p-q-1)/2}\exp\left(\frac{1}{2}\sum_{i=1}^{q}\lambda_i\right)I(0 < \lambda_1, \lambda_2, \ldots, \lambda_q < n).$$

*Then, $e^{-a_n}L_n$ converges weakly to the distribution of $e^{\sigma\xi}$, where $\xi$ is a standard normal, and*

$$a_n = \frac{p^2 q + pq^2}{4n} + \frac{3xy + x^3 y + xy^3}{12} \quad \text{and} \quad \sigma = \frac{xy}{4}.$$

PROOF. Set

(2.11) $$f(x) = \begin{cases} \dfrac{x}{2} + \dfrac{n-p-q-1}{2}\log\left(1 - \dfrac{x}{n}\right), & \text{if } 0 \leq x < n, \\ -\infty, & \text{otherwise.} \end{cases}$$

Then, $L_n = \exp(\sum_{i=1}^{q} f(\lambda_i))$. For any $x \in (0, n)$, by Taylor's expansion, there exists $\xi = \xi_x \in (0, x)$ such that

$$\log\left(1 - \frac{x}{n}\right) = 1 - \frac{x}{n} - \frac{x^2}{2n^2} - \frac{x^3}{3n^3} - \frac{x^4}{4}\cdot\frac{1}{(\xi-n)^4}.$$

Then

(2.12)
$$f(x) = \frac{p+q+1}{2n}x - \frac{n-p-q-1}{4n^2}x^2$$
$$- \frac{n-p-q-1}{6n^3}x^3 + g_n(x)\frac{x^4}{n^3}, \qquad x \in (0, n),$$

where $g_n(x) = -n^3(n - p - q - 1)/(8(\xi - n)^4)$. It is trivial to see that $\sup_{0 \leq x \leq \alpha n} |g_n(x)| \leq (1-\alpha)^{-4}$ for any $\alpha \in (0,1)$. Recall that $\lambda_1, \lambda_2, \ldots, \lambda_q$ are eigenvalues of $X_n' X_n$, where the entries of the $p \times q$ matrix $X_n$ are independent standard normals. Note that $p \sim x\sqrt{n}$ and $q \sim y\sqrt{n}$. By the Theorem from [17] or Theorem 3.1 from [31], there exists a constant $c(x,y) \in (0, \infty)$ such that

(2.13)
$$\frac{\max_{1 \leq i \leq q} \lambda_i}{\sqrt{n}} \to c(x,y)$$

in probability as $n \to \infty$. Define $\Omega_n := \{\max_{1 \leq i \leq q} \lambda_i \leq (c(x,y) + 1)\sqrt{n}\}$. Then

(2.14)
$$P(\Omega_n^c) \to 0$$

as $n \to \infty$. Now on $\Omega_n$, by (2.12),

(2.15)
$$\sum_{i=1}^{q} f(\lambda_i) = \frac{p+q+1}{2n}\operatorname{tr}(X'X) - \frac{n-p-q-1}{4n^2}\operatorname{tr}((X'X)^2)$$
$$- \frac{n-p-q-1}{6n^3}\operatorname{tr}((X'X)^3) + \tilde{g}_n\frac{\operatorname{tr}((X'X)^4)}{n^3},$$

where $|\tilde{g}_n| \in [0, 2)$, as $n$ is sufficiently large. Note that $\operatorname{tr}((X'X)^i)$ are well-defined random variables which do not depend on $\Omega_n$. Easily, $E(\operatorname{tr}(X'X)) = pq$. By Lemma 2.3,

$$E\operatorname{tr}((X'X)^3) \sim pq(p^2 + q^2 + 3pq) \quad \text{and} \quad E\operatorname{tr}((X'X)^4) \leq C(x,y)q^5$$

for some constant $C(x,y)$. It is easy to check that

$$\operatorname{tr}((X'X)^2) = \sum_{j=1}^{q}\sum_{i=1}^{p} x_{ij}^4 + \sum_{j=1}^{q}\sum_{i \neq l=1}^{p} x_{ij}^2 x_{lj}^2$$
$$+ \sum_{i=1}^{p}\sum_{j \neq k=1}^{q} x_{ij}^2 x_{ik}^2 + \sum_{i \neq l, j \neq k} x_{ij}x_{ik}x_{lk}x_{lj}.$$

Then $E\operatorname{tr}((X'X)^2) = pq(p + q + 1)$ [this is sharper than the one corresponding to the case $k = 2$ in (i) of Lemma 2.3]. Now set $h_i = \operatorname{tr}(X'X)^i -$

$E(\mathrm{tr}(X'X)^i)$ for $i = 1, 2, 3$. By simple algebra, we have from (2.15) that

$$\sum_{i=1}^{q} f(\lambda_i) = \frac{p^2 q + pq^2}{4n} + \frac{3xy + x^3 y + xy^3}{12} + O\left(\frac{1}{\sqrt{n}}\right)$$

$$+ \frac{p+q+1}{2n} h_1 - \frac{n-p-q-1}{4n^2} h_2 - \frac{n-p-q-1}{6n^3} h_3 + \tilde{g}_n \frac{h_4}{n^3}$$

on $\Omega_n$ as $n \to \infty$. Recall that $L_n = \exp(\sum_{i=1}^{q} f(\lambda_i))$. By (ii) of Lemma 2.3, both $h_3/n^2$ and $h_4/n^3$ go to zero in probability. By (2.14), to prove the lemma, it suffices to show that

(2.16) $$W_n := \frac{p+q+1}{2n} h_1 - \frac{n-p-q-1}{4n^2} h_2$$

converges to $N(0, \sigma^2)$ weakly,

where $\sigma$ is as in the statement of the lemma. Since $\mathrm{tr}(X'X) = \sum_{i,j} x_{ij}^2$, which is a sum of independent and identically distributed random variables, $\mathrm{Var}(h_1) = \mathrm{Var}(\mathrm{tr}(X'X)) = 2pq$. Therefore, by Lemma 2.4, $\mathrm{Var}(h_2)/n^2$ converges to a positive constant. By Theorem 4.1 from [21], $(h_1/\sqrt{\mathrm{Var}(h_1)}, h_2/\sqrt{\mathrm{Var}(h_2)})$ converges weakly to a normal distribution with mean zero. It follows that $W_n$ converges weakly to a normal distribution with mean zero. We only need to calculate variance $\sigma^2$. Now,

$$\mathrm{Var}(W_n) = \frac{(p+q+1)^2}{4n^2} \mathrm{Var}(\mathrm{tr}(X'X)) + \frac{(n-p-q-1)^2}{16n^4} \mathrm{Var}(\mathrm{tr}((X'X)^2))$$

$$- \frac{(p+q+1)(n-p-q-1)}{4n^3} \cdot \mathrm{Cov}(\mathrm{tr}(X'X), \mathrm{tr}((X'X)^2)).$$

Since $\mathrm{Var}(\mathrm{tr}(X'X)) = 2pq$ as calculated earlier, by Lemma 2.4 again, the above yields

$$\mathrm{Var}(W_n) \to \frac{x^2 y^2}{16},$$

as $n \to \infty$. Therefore, $\sigma^2 = \lim_{n \to \infty} \mathrm{Var}(W_n) = x^2 y^2 / 16$. The proof is completed. $\square$

COROLLARY 2.1. *For $x > 0$ and $y > 0$, let $p_n = [xn^{1/2}]$ and $q_n = [yn^{1/2}]$. Let $f_n(z)$ be the joint probability density function of $Z_n$ as in Theorem 1 and $g_n(z)$ be the joint probability density function of $p_n q_n$ independent standard normals. Then as $n \to \infty$,*

$$\frac{f_n(X_n)}{g_n(X_n)} \text{ converges weakly to } \exp\left(-\frac{x^2 y^2}{8} + \frac{xy}{4}\xi\right),$$

*where $\xi$ and all the entries of $X_n$ are independent standard normals.*

PROOF. Without loss of generality, we assume $y \leq x$. Hence, $q_n \leq p_n$ for any $n \geq 1$. By Lemma 2.5, the density function of $\sqrt{n}Z_n$ is

$$f_n(z) := (\sqrt{2\pi})^{-pq} n^{-pq/2} \frac{\omega(n-p,q)}{\omega(n,q)} \left\{ \det\left(I_q - \frac{z'z}{n}\right)^{(n-p-q-1)/2} \right\} I_0(z'z/n).$$

Obviously, $g_n(z) := (\sqrt{2\pi})^{-pq} e^{-\operatorname{tr}(z'z)/2}$. Let $\lambda_1, \lambda_2, \ldots, \lambda_q$ be the eigenvalues of $X_n' X_n$. Then

$$\frac{f_n(X_n)}{g_n(X_n)} = K_n \cdot L_n,$$

where

$$(2.17) \quad K_n = \left(\frac{2}{n}\right)^{pq/2} \prod_{j=1}^q \frac{\Gamma((n-j+1)/2)}{\Gamma((n-p-j+1)/2)},$$

$$(2.18) \quad L_n = \left\{ \prod_{i=1}^q \left(1 - \frac{\lambda_i}{n}\right) \right\}^{(n-p-q-1)/2} \exp\left(\frac{1}{2} \sum_{i=1}^q \lambda_i\right)$$

if all $\lambda_i$'s are in $(0,n)$, and $L_n$ is zero otherwise. The desired conclusion immediately follows from Lemmas 2.6 and 2.7 on $K_n$ and $L_n$, respectively. □

PROOF OF THEOREM 2. First, we show that the lower bound is strictly between zero and one. Recall $\phi(x,y) = E|\exp(-(x^2y^2/8) + (xy\xi/4)) - 1|$. Then $\phi(x,y) > 0$ because $\xi$ is a nondegenerate random variable. Second, by Hölder's inequality,

$$\phi(x,y) \leq \left\{ E\left[\exp\left(-\frac{x^2y^2}{8} + \frac{xy}{4}\xi\right) - 1\right]^2 \right\}^{1/2}.$$

By expanding the square and using the fact that $E \exp(t\xi) = \exp(t^2/2)$ for any $t \in \mathbb{R}$, we have that

$$\phi(x,y)^2 \leq e^{-x^2y^2/8} - 2e^{-3x^2y^2/32} + 1.$$

Let $\varphi(t) = e^{-t/8} - 2e^{-3t/32} + 1$ for $t \in \mathbb{R}$. Then $\varphi(0) = 0$, $\varphi(+\infty) = 1$ and $\varphi'(t) = (3/16)e^{-t/8}(e^{t/32} - (2/3)) > 0$ for any $t > 0$. Thus, $\phi(x,y) < 1$ for any $x > 0$ and $y > 0$.

Now we prove the remaining part of Theorem 2.

Let us continue to use the notation in Corollary 2.1. First,

$$(2.19) \quad d(\mathcal{L}(\sqrt{n}Z), G_n) = \int_{\mathbb{R}^{pq}} \left|\frac{f_n(z)}{g_n(z)} - 1\right| g_n(z)\,dz = E\left|\frac{f_n(X_n)}{g_n(X_n)} - 1\right|,$$

where $X_n$ has the density function $g_n(z)$, that is, the $pq$ entries of $X_n$ are independent standard normals. Second, by Corollary 2.1,

$$\frac{f_n(X_n)}{g_n(X_n)} \text{ converges weakly to } \exp\left(-\frac{x^2 y^2}{8} + \frac{xy}{4}\xi\right),$$

where $\xi$ is a standard normal. Then, applying Fatou's lemma to (2.19),

$$\liminf_{n \to \infty} d(\mathcal{L}(\sqrt{n}Z), G_n) \geq E\left|\exp\left(-\frac{x^2 y^2}{8} + \frac{xy}{4}\xi\right) - 1\right|.$$

The proof is completed. □

PROOF OF THEOREM 1. Let $p'_n = q'_n = p_n + q_n + [n^{1/4}]$. For an $n$ by $n$ random orthogonal matrix $U$ which has the normalized Haar measure, let $Z_{p,q}$ denote the upper left $p$ by $q$ block of $U$, $1 \leq p$, $q \leq n$. Thus, $Z_{p_n, q_n}$ is a sub-block of $Z_{p'_n, q'_n}$. As a consequence, the joint density function of entries of $Z_{p_n, q_n}$ is a marginal density function of that of $Z_{p'_n, q'_n}$. Therefore, by formula (1.2),

$$(2.20) \qquad \|\mathcal{L}(\sqrt{n}Z_{p_n, q_n}) - G_{p_n q_n}\| \leq \|\mathcal{L}(\sqrt{n}Z_{p'_n, q'_n}) - G_{p'_n q'_n}\|,$$

where $G_{pq}$ is the joint distribution of $pq$ standard normal distributions [one can verify this by choosing $B = A \times \mathbb{R}^{p'_n q'_n - p_n q_n}$ for any Borel set $A \in \mathbb{R}^{p_n q_n}$ and then plugging them into definition (1.2)].

So, to prove the theorem, without loss of generality, we assume $p_n = q_n$ for all $n \geq 1$, $p_n \to \infty$ and $p_n = o(\sqrt{n})$.

As in the proof of Theorem 2,

$$\|\mathcal{L}(\sqrt{n}Z_{p_n, q_n}) - G_{p_n q_n}\| = E|K_n \cdot L_n - 1|,$$

where $K_n$ and $L_n$ are as in (2.17) and (2.18). By following the proof of Lemma 2.6 step by step, we obtain that

$$(2.21) \qquad K_n = \exp\left\{-\frac{p^2 q + pq^2}{4n} + o(1)\right\}$$

as $n \to \infty$. We claim that

$$(2.22) \qquad e^{-(p^2 q + pq^2)/4n} L_n \to 1,$$

in probability as $n \to \infty$. If this is true, then $K_n \cdot L_n \to 1$ in probability as $n \to \infty$. Note that $K_n \cdot L_n \geq 0$ and it is easy to see that $E(K_n \cdot L_n) = \int_{\mathbb{R}^{pq}} f_n(x) \, dx = 1$. These three facts imply that $\{K_n \cdot L_n\}$ is uniformly integrable, that is, $\limsup_{t \to +\infty} \limsup_{n \to \infty} E(K_n L_n I_{\{K_n L_n \geq t\}}) = 0$. It follows that $E|K_n L_n - 1| \to 0$ as $n \to \infty$. The proof is then complete.

Now we prove claim (2.22). Let us go back to the proof of Lemma 2.7. Since $p_n = q_n = o(\sqrt{n})$, the term $c(x, y)$ in (2.13) is equal to zero. So, correspondingly, $\Omega_n = \{\max_{1 \leq i \leq q} \lambda_i \leq \sqrt{n}\}$ and $P(\Omega_n^c) \to 0$ as $n \to \infty$. Recall the definition of $f(x)$ in (2.11) and $L_n = \exp(\sum_{i=1}^q f(\lambda_i))$. On $\Omega_n$, similar to (2.15),

$$\sum_{i=1}^q f(\lambda_i) = \frac{p+q+1}{2n} \operatorname{tr}(X'X) - \frac{n-p-q-1}{4n^2} \operatorname{tr}((X'X)^2)$$

$$+ \tilde{g}_n \frac{\operatorname{tr}((X'X)^3)}{n^2}$$

(2.23)

$$= \frac{p^2q + pq^2}{4n} + \frac{p+q+1}{2n} \cdot h_1 - \frac{n-p-q-1}{4n^2} \cdot h_2$$

$$+ \tilde{g}_n \frac{\operatorname{tr}((X'X)^3)}{n^2},$$

where $\tilde{g}_n$ is a random variable satisfying $|\tilde{g}_n| \in [0, 2)$, as $n$ is sufficiently large, and $h_i = \operatorname{tr}(X'X)^i - E(\operatorname{tr}(X'X)^i)$. Obviously, $h_i$ is well defined on the same probability space as those of $x_{ij}$'s which do not depend on $\Omega_n$. Note that $p = q = o(\sqrt{n})$. Then

(2.24) $$\frac{p}{n} h_1 = \frac{p^2}{n} \cdot \frac{\sum_{i=1}^p \sum_{i=1}^q (x_{ij}^2 - 1)}{p} \to 0,$$

in probability as $n \to \infty$ by the classical central limit theorem of independent and identically distributed random variables. We will show next that the third term on the right-hand side of (2.23) also goes to zero in probability. Indeed,

$$P\left(\frac{|h_2|}{n} \geq \varepsilon\right) \leq \frac{\operatorname{Var}(\operatorname{tr}((X'X)^2))}{n^2 \varepsilon^2} = O\left(\frac{(pq)^2 + 8pq(p+q)^2}{n^2}\right) \to 0$$

by (i) of Lemma 2.4. This says that

(2.25) $$\frac{n-p-q-1}{4n^2} \cdot h_2 \to 0,$$

in probability as $n \to \infty$. Last,

$$\frac{\operatorname{tr}((X'X)^3)}{n^2} = \frac{p^4}{n^2} \cdot \frac{\operatorname{tr}((X'X)^3) - E(\operatorname{tr}((X'X)^3))}{p^4} + \frac{E(\operatorname{tr}((X'X)^3))}{n^2}.$$

By (ii) of Lemma 2.3, the first term on the right-hand side goes to zero in probability. By (i) of Lemma 2.3, $E\operatorname{tr}((X'X)^3) \sim pq(p^2 + q^2 + 3pq)$ as $n \to \infty$. So the second term on the right-hand side goes to zero. Consequently,

(2.26) $$\frac{\operatorname{tr}((X'X)^3)}{n^2} \to 0$$

in probability. Combining (2.23)–(2.26), we obtain

$$\sum_{i=1}^q f(\lambda_i) - \frac{p^2q + pq^2}{4n} \to 0$$

in probability, which, together with the fact that $P(\Omega_n^c) \to 0$, implies (2.22). □

**3. The proof of Theorem 3.** The main tool of proving Theorem 3 is the Gram–Schmidt algorithm. Let us briefly review it first.

Suppose $\{\mathbf{y}_1, \mathbf{y}_2, \ldots, \mathbf{y}_n\}$ is a sequence of $n \times 1$ vectors. Set $\mathbf{w}_1 = \mathbf{y}_1$ and

$$\mathbf{w}_j = \mathbf{y}_j - \sum_{i=1}^{j-1} \frac{\mathbf{y}_j^T \mathbf{w}_i}{\|\mathbf{w}_i\|^2} \mathbf{w}_i, \qquad j = 2, 3, \ldots, n, \tag{3.1}$$

where $\|\mathbf{w}_j\|^2 = \mathbf{w}_j^T \mathbf{w}_j$ $(j = 1, 2, \ldots, n)$. Then, $\{\mathbf{w}_j, 1 \leq j \leq n\}$ are orthogonal, that is, $\mathbf{w}_i^T \mathbf{w}_j = 0$ for any $1 \leq i < j \leq n$. Let $\boldsymbol{\gamma}_j = (1/\|\mathbf{w}_j\|) \mathbf{w}_j$, $j = 1, 2, \ldots, n$. Then the matrix $\boldsymbol{\Gamma}_n = (\boldsymbol{\gamma}_1, \boldsymbol{\gamma}_2, \ldots, \boldsymbol{\gamma}_n)$ is orthonormal. So (3.1) can be rewritten as follows:

$$\mathbf{w}_j = \mathbf{y}_j - \sum_{i=1}^{j-1} (\mathbf{y}_j^T \boldsymbol{\gamma}_i) \boldsymbol{\gamma}_i, \qquad j = 2, 3, \ldots, n. \tag{3.2}$$

The reader is referred to Section A.5 on page 603 from [1] and page 15 from [18] for further details.

Define

$$\boldsymbol{\Delta}_1 = \mathbf{0}, \qquad \boldsymbol{\Delta}_j = \sum_{i=1}^{j-1} (\mathbf{y}_j^T \boldsymbol{\gamma}_i) \boldsymbol{\gamma}_i \quad \text{and} $$

$$L_j = \left| \sqrt{\frac{n}{\|\mathbf{w}_j\|^2}} - 1 \right|, \qquad j = 1, 2, \ldots, n. \tag{3.3}$$

Note $\mathbf{y}_j^T \boldsymbol{\gamma}_i \in \mathbb{R}^1$ and rewrite $(\mathbf{y}_j^T \boldsymbol{\gamma}_i) \boldsymbol{\gamma}_i = (\boldsymbol{\gamma}_i \boldsymbol{\gamma}_i^T) \mathbf{y}_j$. It is easy to check that

$$\mathbf{w}_j = (\mathbf{I}_n - \boldsymbol{\Gamma}_{n,j} \boldsymbol{\Gamma}_{n,j}^T) \mathbf{y}_j, \qquad \boldsymbol{\Delta}_j = \boldsymbol{\Gamma}_{n,j} \boldsymbol{\Gamma}_{n,j}^T \mathbf{y}_j \quad \text{and}$$

$$\boldsymbol{\gamma}_j = \frac{\mathbf{y}_j}{\sqrt{n}} - \frac{\boldsymbol{\Delta}_j}{\sqrt{n}} + \mathbf{u}_j, \tag{3.4}$$

where $\boldsymbol{\Gamma}_{n,j} = (\boldsymbol{\gamma}_1, \boldsymbol{\gamma}_2, \ldots, \boldsymbol{\gamma}_{j-1})$ and $\mathbf{u}_j = (1 - n^{-1/2} \|\mathbf{w}_j\|) \boldsymbol{\gamma}_j$.

One repeatedly used fact in later proofs is that if the $n^2$ elements of $\mathbf{Y} = (\mathbf{y}_1, \mathbf{y}_2, \ldots, \mathbf{y}_n)$ are i.i.d. standard normals, then $\boldsymbol{\Gamma}_n = (\boldsymbol{\gamma}_1, \boldsymbol{\gamma}_2, \ldots, \boldsymbol{\gamma}_n)$ follows the normalized Haar measure on the orthogonal group $O(n)$. In particular, $\boldsymbol{\gamma}_i$'s are identically distributed and

$$\mathcal{L}(\boldsymbol{\gamma}_i) = \mathcal{L}\left( \frac{\mathbf{y}_1}{\|\mathbf{y}_1\|} \right) \tag{3.5}$$

for any $i = 1, 2, \ldots, n$.

For any $n \times n$ orthogonal matrix $\mathbf{G}$, observe that $\mathcal{L}(\mathbf{G}\mathbf{\Gamma}_n^{-1}) = \mathcal{L}((\mathbf{\Gamma}_n\mathbf{G}^T)^{-1}) = \mathcal{L}(\mathbf{\Gamma}_n^{-1})$ by the invariance property of Haar measures. Also, $\mathbf{\Gamma}_n^{-1} = \mathbf{\Gamma}_n^T$. From the uniqueness of Haar measures, we obtain another useful fact that

$$\mathcal{L}(\mathbf{\Gamma}_n) = \mathcal{L}(\mathbf{\Gamma}_n^T). \tag{3.6}$$

We will use the following notation. Let $A = (a_{ij})$ be a $p$ by $q$ matrix. Then

$$\|A\| := \max_{1 \leq i \leq p, 1 \leq j \leq q} |a_{ij}|. \tag{3.7}$$

The following definition will also be used:

$$\varepsilon_n(m) = \max_{1 \leq i \leq n, 1 \leq j \leq m} |\sqrt{n}\gamma_{ij} - y_{ij}| \quad \text{and}$$

$$n_\alpha = \left[\frac{n}{\log n - (5/4)\log(\log n)}\alpha\right] \tag{3.8}$$

for $\alpha > 0$ and $n \geq 2$.

The following says that, to prove part (iii) of Theorem 3, we only need to work on $\max_{2 \leq j \leq m} \|\mathbf{\Delta}_j\|$.

LEMMA 3.1. *Let $\varepsilon_n(m)$ and $n_\alpha$ be as in (3.8). Then*

$$P\left(\left|\varepsilon_n(n_\alpha) - \max_{2 \leq j \leq n_\alpha} \|\mathbf{\Delta}_j\|\right| \geq \delta\right) \to 0,$$

*as $n \to \infty$ for any $\alpha > 0$ and $\delta > 0$.*

The following lemma is the key in the proof of Theorem 3. A recursive inequality is derived. It implies that all $\mathbf{\Delta}_j$'s are almost independent when $j \leq n_\alpha$.

LEMMA 3.2. *Let $\xi$ be a standard normal. Given $\alpha > 0$ and $t > 0$, define*

$$f_n^+(k) = P\left(|\xi| > t\left(\sqrt{\frac{n}{k}} + \frac{(\log n)^8}{\sqrt{n}}\right)\right), \quad k = 1, 2, \ldots, n,$$

*and $f_n^-(k)$ as the probability above when "+" on the right-hand side is replaced by "−." Then there exists a constant $C = C_{\alpha,t} > 0$ such that $P(\max_{2 \leq j \leq k+1} \|\mathbf{\Delta}_j\| \leq t)$ is bounded below and above, respectively, by*

$$(1 - nf_n^-(k))P\left(\max_{2 \leq j \leq k} \|\mathbf{\Delta}_j\| \leq t\right) - \frac{(\log n)^C}{n^{(t^2/\alpha)-2}}$$

*and*

$$(1 - nf_n^+(k))P\left(\max_{2 \leq j \leq k} \|\mathbf{\Delta}_j\| \leq t\right) + \frac{(\log n)^C}{n^{(t^2/\alpha)-2}},$$

*uniformly on $n/(\log n)^3 \leq k \leq n_\alpha$ as $n$ is sufficiently large, where $n_\alpha$ is as in (3.8).*

PROOF OF THEOREM 3. Part (i) is obvious. As for (ii), take $r = 1/\log n$, $s = (\log n)^{3/4}$, $t = t$, $m = m'_n = [\delta n/\log n]$ for some $\delta < \min\{1/4, t^2/100\}$ in Lemma A.4. Trivially, $t^2/(3(m + \sqrt{n})) \geq t^2(\log n)/(6n\delta)$ and $1/s \leq 1$, as $n$ is sufficiently large. We obtain that

$$P(\varepsilon_n(m'_n) \geq 3t)$$
$$\leq 4ne^{-n/(4\log n)^2} + 3n^2 e^{-(\log n)^{3/2}/2} + \frac{3n^2}{t}\left(1 + \frac{t^2}{6\delta}\frac{\log n}{n}\right)^{-n/2}$$
$$\to 0,$$

as $n \to \infty$ by the choice of $\delta$.

Now we prove (iii). To simplify notation, set $m = n_\alpha$. We actually will show that

$$(3.9) \qquad P\left(\max_{2\leq j \leq m} \|\mathbf{\Delta}_j\| \leq t\right) \to \begin{cases} 1, & \text{if } t > 2\sqrt{\alpha}, \\ e^{-Kt^2}, & \text{if } t = 2\sqrt{\alpha}, \\ 0, & \text{if } t \in (\sqrt{3\alpha}, 2\sqrt{\alpha}), \end{cases}$$

where $K = (8\sqrt{2\pi})^{-1}$. Since $P(\max_{2\leq j \leq m} \|\mathbf{\Delta}_j\| \leq t)$ is increasing in $t$, the above implies that the left-hand side above goes to zero for any $t \in (0, 2\sqrt{\alpha})$. This together with (3.9) implies that $\max_{2\leq j \leq m} \|\mathbf{\Delta}_j\|$ converges to $2\sqrt{\alpha}$ in probability. Lemma 3.1 says that $\varepsilon_n(n_\alpha) - \max_{2\leq j \leq n_\alpha} \|\mathbf{\Delta}_j\|$ converges to zero in probability as $n \to \infty$. It follows that

$$(3.10) \qquad \varepsilon_n(n_\alpha) \to 2\sqrt{\alpha},$$

in probability as $n \to \infty$. We next show that this implies that $\varepsilon_n([n\alpha/\log n]) \to 2\sqrt{\alpha}$ as $n \to \infty$. Indeed, set $k_\alpha = [n\alpha/\log n]$. For any $\delta \in (0, \sqrt{\alpha})$, choose $\alpha_1$ such that

$$\left(\sqrt{\alpha} - \frac{\delta}{4}\right)^2 < \alpha_1 < \alpha.$$

Then $n_{\alpha_1} < k_\alpha \leq n_\alpha$, as $n$ is sufficiently large. It follows from the definition of $\varepsilon_n(m)$ that $\varepsilon_n(n_{\alpha_1}) \leq \varepsilon_n(k_\alpha) \leq \varepsilon_n(n_\alpha)$, as $n$ is sufficiently large. Therefore,

$$P(|\varepsilon_n(k_\alpha) - 2\sqrt{\alpha}| > \delta) \leq P(\varepsilon_n(k_\alpha) > 2\sqrt{\alpha} + \delta) + P(\varepsilon_n(k_\alpha) < 2\sqrt{\alpha} - \delta)$$
$$\leq P(\varepsilon_n(n_\alpha) > 2\sqrt{\alpha} + \delta) + P\left(\varepsilon_n(n_{\alpha_1}) < 2\sqrt{\alpha_1} - \frac{\delta}{2}\right),$$

as $n$ is sufficiently large. The above two terms go to zero as $n \to \infty$ by (3.10). Then (iii) follows.

Now we show (3.9).

We continue to use the notation in Lemma 3.2. Set

$$A_k = P\left(\max_{2\leq j \leq k} \|\mathbf{\Delta}_j\| \leq t\right), \qquad b_k^+ = 1 - nf_n^+(k), \qquad b_k^- = 1 - nf_n^-(k),$$

$$c_n = \frac{(\log n)^C}{n^{(t^2/\alpha)-2}} \quad \text{and} \quad m' = \left[\frac{n}{(\log n)^3}\right] + 2.$$

By Lemma A.1, $P(|\xi| \geq x) \sim (2/(\sqrt{2\pi}x))\exp(-x^2/2)$ as $x \to +\infty$ for a standard normal $\xi$. Here and later, the notation "$f(x) \sim g(x)$ as $x \to +\infty$" means that $\lim_{x \to +\infty} f(x)/g(x) = 1$. The same interpretation applies to $\alpha_n \sim \beta_n$ as $n \to \infty$. It is easy to check that

$$(3.11) \qquad \text{both } f_n^+(k) \text{ and } f_n^-(k) \sim \frac{2}{t\sqrt{2\pi}}\left(\frac{k}{n}\right)^{1/2} e^{-(t^2/2)(n/k)}$$

uniformly on $m' \leq k \leq m$ as $n \to \infty$, and also that

$$(3.12) \qquad 1 > \max\{b_i^+, b_i^-;\ m' \leq i \leq m\} \to 1$$

as $n \to \infty$, provided $t > \sqrt{2\alpha}$. By Lemma 3.2,

$$b_k^- A_{k-1} - c_n \leq A_k \leq b_k^+ A_{k-1} + c_n$$

for all $m' \leq k \leq m$, as $n$ is sufficiently large. By iteration, we obtain

$$(3.13) \qquad A_m \geq \left(\prod_{j=m'}^{m} b_j^-\right) A_{m'-1} - c_n \sum_{j=0}^{m-m'+2} \left[\max_{m' \leq i \leq m}\{b_i^-\}\right]^j.$$

By (3.12), the second term on the right-hand side is no larger than $nc_n \leq (\log n)^C/n^{(t^2/\alpha)-3}$, as $n$ is sufficiently large. Further, applying the same argument in (3.13) to the "+" case, we obtain

$$(3.14) \left(\prod_{j=m'}^{m} b_j^-\right) A_{m'-1} - \frac{(\log n)^C}{n^{(t^2/\alpha)-3}} \leq A_m \leq \left(\prod_{j=m'}^{m} b_j^+\right) A_{m'-1} + \frac{(\log n)^C}{n^{(t^2/\alpha)-3}},$$

as $n$ is sufficiently large. By definition, $A_k = P(\max_{2 \leq j \leq k}\|\mathbf{\Delta}_j\| \leq t)$. From the proved (ii), we know that $A_{m'-1} \to 1$ as $n \to \infty$ for any $t > 0$. Evidently, $(\log n)^C n^{3-(t^2/\alpha)} \to 0$, provided $t > \sqrt{3\alpha}$. So to prove (3.9), we only need to show that

$$(3.15) \text{ both } \prod_{j=m'}^{m} b_j^- \text{ and } \prod_{j=m'}^{m} b_j^+ \to \begin{cases} 1, & \text{if } t > 2\sqrt{\alpha}, \\ e^{-Kt^2}, & \text{if } t = 2\sqrt{\alpha}, \\ 0, & \text{if } t \in (\sqrt{3\alpha}, 2\sqrt{\alpha}), \end{cases}$$

as $n \to \infty$. Recall $b_j^+ = 1 - nf_n^+(k)$ and $b_j^- = 1 - nf_n^-(k)$. Since $|\log(1+x) - x| \leq x^2$ for $x$ small enough, by (3.11) and (3.12),

$$\prod_{j=m'}^{m} b_j^+ \leq \exp\left(-n \sum_{k=m'}^{m} f_n^+(k)\right) \cdot \exp\left(+n^2 \sum_{k=m'}^{m} f_n^+(k)^2\right),$$

$$\prod_{j=m'}^{m} b_j^- \geq \exp\left(-n \sum_{k=m'}^{m} f_n^-(k)\right) \cdot \exp\left(-n^2 \sum_{k=m'}^{m} f_n^-(k)^2\right),$$

as $n$ is sufficiently large. Also, the fact $f_n^+(k) \leq f_n^-(k)$ implies that $b_j^+ \geq b_j^-$. So (3.15) is reduced to show that

$$n^2 \sum_{k=m'}^{m} f_n^+(k)^2 \to 0 \quad \text{and}$$

(3.16)
$$n \sum_{k=m'}^{m} f_n^+(k) \to \begin{cases} 0, & \text{if } t > 2\sqrt{\alpha}, \\ Kt^2, & \text{if } t = 2\sqrt{\alpha}, \\ +\infty, & \text{if } t \in (\sqrt{3\alpha}, 2\sqrt{\alpha}), \end{cases}$$

and that the above is also true if $f_n^+(k)$ is replaced by $f_n^-(k)$.

By (3.11) again, $n^2 \sum_{k=m'}^{m} f_n^+(k)^2 \leq (\log n)^C n^{3-(t^2/\alpha)} \to 0$ as $n \to \infty$, provided $t > \sqrt{3\alpha}$. Similarly, $n^2 \sum_{k=m'}^{m} f_n^-(k)^2 \to 0$ for $t > \sqrt{3\alpha}$. Let

$$g(x) = \frac{2}{t\sqrt{2\pi}} x^{1/2} e^{-t^2/(2x)}$$

for $x > 0$. By the uniform convergence of $f_n^+(k)/g(k/n)$ and $f_n^-(k)/g(k/n)$ as $n \to \infty$ over $k \in [m', m]$ as in (3.11), to prove the second part in (3.16), it is enough to show

(3.17) $$n \sum_{k=m'}^{m} g\left(\frac{k}{n}\right) \text{ goes to the second limit in (3.16)}$$

as $n \to \infty$. Note that $g(x)$ is nonnegative and increasing in $x$ over $[0, +\infty)$, it is elementary to see that

$$\left| \frac{1}{n} \sum_{k=m'}^{m} g\left(\frac{k}{n}\right) - \int_0^{m/n} g(x)\, dx \right| \leq \int_{m/n}^{(m+1)/n} g(x)\, dx + \int_0^{m'/n} g(x)\, dx.$$

Using $\sqrt{x} e^{-t^2/(2x)} \leq e^{-t^2/(2x)}$ on $x \in [0,1]$, the first integral on the right-hand side is bounded by $(1/n) \exp(-nt^2/(2m+2)) \leq n^{-1-(t^2/(2\alpha))} (\log n)^C$, as $n$ is sufficiently large; the second one is bounded by $\exp(-(\log n)^2)$ as $n$ is large because $m' \sim n(\log n)^{-3}$ by definition. Hence,

(3.18) $$\frac{1}{n} \sum_{k=m'}^{m} g\left(\frac{k}{n}\right) - \int_0^{m/n} g(x)\, dx = o\left(\frac{1}{n^2}\right)$$

as $n \to \infty$ if $t > \sqrt{2\alpha}$. Now we evaluate the integral.

Write $\sqrt{x} \exp(-t^2/(2x))\, dx = (2t^{-2} x^{5/2}) d(e^{-t^2/(2x)})$. By integration by parts,

$$I_n := \int_0^{m/n} \sqrt{x} e^{-t^2/(2x)}\, dx = \frac{2}{t^2} \left(\frac{m}{n}\right)^{5/2} e^{-nt^2/(2m)} - \frac{5}{t^2} \int_0^{m/n} \sqrt{x^3} e^{-t^2/(2x)}\, dx.$$

Note that $\sqrt{x^3} \leq (m/n)\sqrt{x}$ on $[0, m/n]$. The last integral is less than or equal to $(m/n)I_n$. But $m/n \to 0$, thus,

$$I_n \sim \frac{2}{t^2} \left(\frac{m}{n}\right)^{5/2} e^{-nt^2/(2m)}.$$

By the definition of $m$, $nt^2/(2m) = (t^2/(2\alpha))(\log n - (5/4)\log_2 n) + O(n^{-1}(\log n)^2)$ as $n \to \infty$. It follows that

$$(3.19) \qquad I_n \sim \frac{2\alpha^{5/2}}{t^2} \cdot \frac{1}{n^{t^2/(2\alpha)}} (\log n)^{5t^2/(8\alpha)-5/2}.$$

From (3.18),

$$n \sum_{k=m'}^{m} g\left(\frac{k}{n}\right) \sim n^2 \int_0^{m/n} g(x)\,dx$$

$$= \frac{2n^2}{\sqrt{2\pi}t} \cdot I_n \sim \frac{4\alpha^{5/2}}{\sqrt{2\pi}t^3} \cdot \frac{1}{n^{t^2/(2\alpha)-2}} \cdot (\log n)^{5t^2/(8\alpha)-5/2},$$

provided $t > \sqrt{3\alpha}$. Recall $K = (8\sqrt{2\pi})^{-1}$. The above implies (3.17). □

**4. Technical lemmas.** Now we prove the lemmas used in the previous sections. To see them clearly, we break them into two subsections.

4.1. *The proofs of lemmas used in Section* 2.

PROOF OF LEMMA 2.1.  (i) First, when $n = 1$, $\Gamma(n+(1/2))/(\sqrt{n}\Gamma(n)) = \sqrt{\pi}/2 \in (5/6, 1)$. So (i) is true for $n = 1$. Now assume $n \geq 2$.

Using the fact that $\Gamma(x+1) = x\Gamma(x)$ for any $x > 0$ and $\Gamma(1/2) = \sqrt{\pi}$, we have that

$$\frac{\Gamma(n+(1/2))}{\Gamma(n)} = \frac{\sqrt{\pi}n}{2^{2n}} \cdot \frac{(2n)!}{(n!)^2}.$$

By Stirling's formula (see, e.g., Lemma 1 on page 45 from [6]), $n! = \sqrt{2\pi n}n^n \cdot e^{-n+\theta_n/(12n)}$ for all $n \geq 2$, where

$$(4.1) \qquad \frac{n}{n+1/12} < \theta_n < 1.$$

It is easily checked that

$$(4.2) \qquad \frac{\Gamma(n+(1/2))}{\sqrt{n}\,\Gamma(n)} = \exp\left(\frac{\theta_n - 4\theta'_n}{24n}\right)$$

for some $\theta_n$ corresponding to $2n$ and $\theta'_n$ corresponding to $n$ in (4.1). Evidently, $(\theta_n - 4\theta'_n)/24 \in (-1/6, 0)$ for all $n \geq 2$. Then the desired result follows by using the inequality $e^x > 1 + x$ for all $x \neq 0$.

(ii) A direct verification shows that (ii) is true for $n = 1$. Now assume $n \geq 2$. If $n = 2k$ for some integer $k \geq 1$, then (ii) follows from (i). Now suppose $n = 2k+1$ for $k \geq 1$. Trivially,

$$\frac{\Gamma((n+1)/2)}{\sqrt{n/2}\,\Gamma(n/2)} = \left(\frac{\Gamma(k+(1/2))}{\sqrt{k}\,\Gamma(k)}\right)^{-1} \cdot \sqrt{\frac{2k}{2k+1}}.$$

By (i), the above ratio is between $\sqrt{2k/(2k+1)}$ and $(1-(6k)^{-1})^{-1}$. By a simple calculation, $\sqrt{2k/(2k+1)} \geq 1 - (3/5n)$ and $(1-(6k)^{-1})^{-1} \leq 1 + (5k)^{-1}$ for all $k \geq 1$. So (ii) follows. $\square$

PROOF OF LEMMA 2.2. By the multivariate Taylor's expansion formula (see page 361 from [2] and page 172 from [22]),

$$f(x,y) = f\left(\frac{i}{n}, \frac{j}{n}\right) + f'_x\left(\frac{i}{n}, \frac{j}{n}\right)\left(x - \frac{i}{n}\right) + f'_y\left(\frac{i}{n}, \frac{j}{n}\right)\left(y - \frac{j}{n}\right) + \delta_{ij}(\xi, \eta),$$

for some $\xi \in [i/n, x]$ and $\eta \in [j/n, y]$, where

$$(4.3) \quad \delta_{ij}(x,y) = \frac{1}{2}\left(\left(x - \frac{i}{n}\right)^2 \frac{\partial^2 f}{\partial x^2} + 2\left(x - \frac{i}{n}\right)\left(y - \frac{j}{n}\right)\frac{\partial^2 f}{\partial x \partial y} + \left(y - \frac{j}{n}\right)^2 \frac{\partial^2 f}{\partial y^2}\right).$$

By the given condition,

$$|\delta_{ij}(x,y)| \leq \frac{M}{2}\left(\left(x - \frac{i}{n}\right) + \left(y - \frac{j}{n}\right)\right)^2 \leq M\left(\left(x - \frac{i}{n}\right)^2 + \left(y - \frac{j}{n}\right)^2\right).$$

Then

$$\int_{j/n}^{(j+1)/n} \int_{i/n}^{(i+1)/n} f(x,y)\,dx\,dy$$
$$= \frac{1}{n^2} f\left(\frac{i}{n}, \frac{j}{n}\right) + \frac{1}{2n^3} f'_x\left(\frac{i}{n}, \frac{j}{n}\right) + \frac{1}{2n^3} f'_y\left(\frac{i}{n}, \frac{j}{n}\right) + \delta'_{ij},$$

where

$$|\delta'_{ij}| = \left|\int_{j/n}^{(j+1)/n}\int_{i/n}^{(i+1)/n} \delta_{ij}(\xi,\eta)\,dx\,dy\right| \leq M\int_0^{1/n}\int_0^{1/n}(x^2+y^2)\,dx\,dy = \frac{2M}{3n^4}$$

since $|\delta_{ij}(\xi,\eta)| \leq M((x-i/n)^2 + (y-j/n)^2)$ by (4.3). The desired result follows by taking the sum over $i$ from $i_1$ to $i_2$, and $j$ from $j_1$ to $j_2$. $\square$

PROOF OF LEMMA 2.4. (i) It is not difficult to check that

$$\operatorname{tr}(X'X) = \sum_{j=1}^{q}\sum_{i=1}^{p} x_{ij}^2;$$

$$(4.4) \quad \operatorname{tr}((X'X)^2) = \sum_{j=1}^{q}\sum_{i=1}^{p} x_{ij}^4 + \sum_{j=1}^{q}\sum_{i\neq l=1}^{p} x_{ij}^2 x_{lj}^2$$
$$+ \sum_{i=1}^{p}\sum_{j\neq k=1}^{q} x_{ij}^2 x_{ik}^2 + \sum_{i\neq l, j\neq k} x_{ij}x_{ik}x_{lk}x_{lj}.$$

Let

$$B_1 = \sum_{j=1}^{q}\sum_{i=1}^{p}(x_{ij}^4 - 3), \qquad B_2 = \sum_{j=1}^{q}\sum_{i\neq l=1}^{p}(x_{ij}^2-1)(x_{lj}^2-1),$$

$$B_3 = \sum_{i=1}^{p}\sum_{j\neq k=1}^{q}(x_{ij}^2-1)(x_{ik}^2-1), \qquad B_4 = \sum_{i\neq l, j\neq k} x_{ij}x_{ik}x_{lk}x_{lj}.$$

By a simple algebra,

$$(4.5) \qquad \operatorname{tr}((X'X)^2) = \left(\sum_{i=1}^{4} B_i\right) + 2(p+q-2)\operatorname{tr}(X'X) + C_{p,q},$$

where $C_{p,q}$ is a constant on $p$ and $q$. It is easy to check that $EB_i = 0$ for $1 \le i \le 4$, $\operatorname{Cov}(B_i, B_j) = 0$ for all $1 \le i \neq j \le 4$, and $\operatorname{Cov}(B_i, \operatorname{tr}(X'X)) = 0$ for $i = 2, 3, 4$. Also, each $B_i$ is a sum of uncorrelated random variables. Therefore,

$$\operatorname{Var}(\operatorname{tr}((X'X)^2)) = \left(\sum_{i=1}^{4}\operatorname{Var}(B_i)\right) + 4(p+q-2)^2\operatorname{Var}(\operatorname{tr}(X'X))$$
$$+ 2\operatorname{Cov}(B_1, \operatorname{tr}(X'X)).$$

Now it is easy to verify that $\operatorname{Cov}(B_1, \operatorname{tr}(X'X)) = O(p^2)$ and $\operatorname{Var}(B_i) = O(p^3)$ for $i = 1, 2, 3$ as $p \to \infty$. Moreover, $\operatorname{Var}(B_4) = pq(p-1)(q-1)$ and $\operatorname{Var}(\operatorname{tr}(X'X)) = 2pq$. Combining these quantities together, we obtain (i).

(ii) By (4.5) again,

$$\operatorname{Cov}(\operatorname{tr}(X'X), \operatorname{tr}((X'X)^2)) = \operatorname{Cov}(\operatorname{tr}(X'X), B_1) + 2(p+q-2)\cdot \operatorname{Var}(\operatorname{tr}(X'X))$$
$$\sim 4pq(p+q)$$

as $n \to \infty$. □

4.2. *The proofs of lemmas used in Section* 3. Before the proof of these lemmas, we need some preliminary results for a preparation.

LEMMA 4.1. *Let $E_i$, $i = 0, 1, 2, \ldots, n$, be events in a probability space $(\Omega, \mathcal{F}, P)$. Then*

$$\left|P\left(\bigcap_{i=0}^{n} E_i\right) - P(E_0) + \sum_{i=1}^{n} P(E_0 \backslash E_i)\right| \le \sum_{1 \le i < j \le n} P(E_i^c E_j^c).$$

PROOF. First, $P(E_0) - P(\bigcap_{i=0}^{n} E_i) = P(\bigcup_{i=1}^{n} E_0 \backslash E_i)$. By Bonferoni's inequality, it is bounded above and below, respectively, by

$$\sum_{i=1}^{n} P(E_0 \backslash E_i) \quad \text{and} \quad \sum_{i=1}^{n} P(E_0 \backslash E_i) - \sum_{1 \le i < j \le n} P((E_0 \backslash E_i) \cap (E_0 \backslash E_j)).$$

Note that $(E_0\backslash E_i)\cap (E_0\backslash E_j)\subset E_i^c E_j^c$. The desired conclusion follows. $\square$

LEMMA 4.2. *Let $\{\xi_i; i\geq 1\}$ be a sequence of i.i.d. random variables with the standard normal distribution. Set $S_k = \sum_{i=1}^k \xi_i^2$. Then*

$$P\left(\left|\frac{S_n}{S_m} - \frac{n}{m}\right| \geq x\right) \leq 6\exp\left(-\frac{m^4 x^2}{48n^3}\right)$$

*for any $m\geq 1$, $n\geq 1$ and $x > 0$ satisfying $m\leq n/2$ and $x\leq n/m$.*

PROOF. Write

$$\frac{S_n}{S_m} - \frac{n}{m} = \frac{(m-n)(S_m - m) + m[(S_n - S_m) - (n-m)]}{mS_m}.$$

Then

$$\left|\frac{S_n}{S_m} - \frac{n}{m}\right| \leq \frac{n}{mS_m} \max\{|S_m - m|,\ |(S_n - S_m) - (n-m)|\}.$$

Since the distribution of $S_n - S_m$ is equal to that of $S_{n-m}$, we have that

$$P\left(\left|\frac{S_n}{S_m} - \frac{n}{m}\right| \geq x\right)$$

(4.6)
$$\leq P\left(S_m \leq \frac{m}{2}\right) + P\left(|S_m - m| > \frac{m^2 x}{2n}\right)$$
$$+ P\left(|S_{n-m} - (n-m)| > \frac{m^2 x}{2n}\right).$$

Let $P_1, P_2$ and $P_3$ stand for the previous three probabilities in order. Define $I(x) := \sup_{\theta\in\mathbb{R}}\{\theta x - \log(E\exp(\theta\xi_1^2))\}$ for $x\in\mathbb{R}$. It is not difficult to verify the following:

(i) $I(x) = (x - 1 - \log x)/2$ for $x > 0$; $I(x) = +\infty$ for $x \leq 0$;
(ii) $I(x)$ is increasing on $[1,\infty)$ and decreasing on $(0,1)$.

The above two facts can be also seen in Lemma 3.2 from [19]. By (i) of Lemma A.3,

$$P_1 \leq 2e^{-mI(1/2)} \leq 2\exp(-(\log 4 - 1)m/4) \leq 2\exp(-m/12).$$

Define $\eta(x) = x - \log(1 + x) - (x^2/3)$ for $x > -1$. Then $\eta(0) = 0$ and $\eta'(x) = x(1 - 2x)(1 + x)^{-1}/3$. Hence, $\eta'(x) \geq 0$ for $x\in[0, 1/2]$ and $\eta'(x) < 0$ for $x\in[-1/2, 0)$. It follows that $x - \log(1+x) \geq x^2/3$ for $|x| < 1/2$. Therefore,

$$P_2 \leq 2\exp\left\{-m\cdot \max\left\{I\left(1 + \frac{mx}{2n}\right),\ I\left(1 - \frac{mx}{2n}\right)\right\}\right\}$$
$$\leq 2e^{-m^3 x^2/(24n^2)},$$

provided $x \leq n/m$, where property (ii) of $I(x)$ above is used. Similarly,

$$P_3 \leq P\left(\left|\frac{S_{n-m}}{n-m} - 1\right| > \frac{m^2 x}{2n^2}\right)$$
$$\leq 2\exp\left\{-(n-m) \cdot \max\left\{I\left(1 + \frac{m^2 x}{2n^2}\right), I\left(1 - \frac{m^2 x}{2n^2}\right)\right\}\right\}$$
$$\leq 2e^{-m^4 x^2/(48n^3)},$$

provided $m \leq n/2$ and $x \leq n^2/m^2$, where the fact that $n - m \geq n/2$ is used in the last step. Thus,

$$P_1 + P_2 + P_3 \leq 6\exp\left(-\min\left\{\frac{m}{12}, \frac{m^3 x^2}{24n^2}, \frac{m^4 x^2}{48n^3}\right\}\right)$$

if $m \leq n/2$ and $x \leq n/m$. By a simple verification, the minimum above is actually $m^4 x^2/(48n^3)$. This together with (4.6) proves the lemma. □

PROOF OF LEMMA 3.1. Write $m = n_\alpha$ for simplification. By (3.4), we know that

$$\max_{1 \leq j \leq m} \|\!|\sqrt{n}\boldsymbol{\gamma}_j - \mathbf{y}_j + \boldsymbol{\Delta}_j|\!\| \leq \max_{1 \leq j \leq m} \|\!|\sqrt{n}\mathbf{u}_j|\!\|,$$

where $\mathbf{u}_j = (1 - n^{-1/2}\|\mathbf{w}_j\|)\boldsymbol{\gamma}_j$. By the triangle inequality,

$$\left|\varepsilon_n(m) - \max_{2 \leq j \leq m} \|\!|\boldsymbol{\Delta}_j|\!\|\right| \leq \max_{1 \leq j \leq m} \|\!|\sqrt{n}\mathbf{u}_j|\!\|$$
$$\leq \left\{\max_{1 \leq j \leq n} \|\!|\sqrt{n}\boldsymbol{\gamma}_j|\!\|\right\} \cdot \max_{1 \leq j \leq m}\left|1 - \frac{\|\mathbf{w}_j\|^2}{n}\right|,$$

where the inequality $|1 - \sqrt{x}| \leq |1 - x|$ is used in the last step. Proposition 1 from [19] implies that

$$\sqrt{\frac{n}{\log n}} \max_{1 \leq j \leq n} \|\!|\boldsymbol{\gamma}_j|\!\| \overset{P}{\to} 2$$

as $n \to \infty$. To prove the lemma, it suffices to show that

(4.7) $$B_n := \sqrt{\log n} \max_{1 \leq j \leq m}\left|1 - \frac{\|\mathbf{w}_j\|^2}{n}\right| \to 0,$$

in probability as $n \to \infty$. By orthogonality, $(I_n - \boldsymbol{\Gamma}_{n,j}\boldsymbol{\Gamma}_{n,j}^T)^2 = I_n - \boldsymbol{\Gamma}_{n,j}\boldsymbol{\Gamma}_{n,j}^T$. This says that $I_n - \boldsymbol{\Gamma}_{n,j}\boldsymbol{\Gamma}_{n,j}^T$ is an idempotent matrix. So by (3.4), $\mathbf{w}_j \sim N_n(\mathbf{0}, I_n - \boldsymbol{\Gamma}_{n,j}\boldsymbol{\Gamma}_{n,j}^T)$ conditioning on $\mathbf{y}_1, \mathbf{y}_2, \ldots, \mathbf{y}_{j-1}$, where $\boldsymbol{\Gamma}_{n,j} = (\boldsymbol{\gamma}_1, \boldsymbol{\gamma}_2, \ldots, \boldsymbol{\gamma}_{j-1})$. In this context, "∼" means that both sides of "∼" have the same probability distribution. It also follows that $\text{rank}(I_n - \boldsymbol{\Gamma}_{n,j}\boldsymbol{\Gamma}_{n,j}^T) = \text{trace}(I_n -$ $\mathbf{\Gamma}_{n,j}\mathbf{\Gamma}_{n,j}^T) = \text{trace}(I_n) - \text{trace}(\mathbf{\Gamma}_{n,j}\mathbf{\Gamma}_{n,j}^T) = n - j + 1$. By Lemma A.2, $\|\mathbf{w}_j\|^2 \sim \chi^2(n - j + 1)$. Obviously, $2tn/\sqrt{\log n} - j \geq tn/\sqrt{\log n}$ for all $1 \leq j \leq m$, as $n$ is sufficiently large. Let $\{\xi_1, \xi_2, \ldots, \xi_n\}$ be independent standard normals. Then

$$P\left(\left|1 - \frac{\|\mathbf{w}_j\|^2}{n}\right| \geq 2t(\log n)^{-1/2}\right)$$

(4.8)
$$\leq P\left(\left|\sum_{k=1}^{n-j+1}(\xi_k^2 - 1)\right| \geq \frac{tn}{\sqrt{\log n}}\right)$$

$$\leq P\left(\frac{1}{(n-j+1)^{1/2}}\left|\sum_{k=1}^{n-j+1}(\xi_k^2 - 1)\right| \geq n^{1/3}\right)$$

$$\leq \exp(-n^{1/3}),$$

uniformly for $1 \leq j \leq m$ as $n$ is sufficiently large, where Lemma A.3 is used in the last inequality [heuristically, since $\sum_{k=1}^{n-j+1}(\xi_k^2 - 1)$ is a sum of i.i.d. random variables with mean zero and variance equal to two, one can think of $\sum_{k=1}^{n-j+1}(\xi_k^2 - 1)/\sqrt{n-j+1}$ as a normal. Then the last inequality above is intuitive]. By the union bound,

$$P(B_n \geq 2t) \leq n \cdot \max_{1 \leq j \leq m} P\left(\left|1 - \frac{\|\mathbf{w}_j\|^2}{n}\right| \geq 2t(\log n)^{-1/2}\right) \leq n \cdot \exp(-n^{1/3}) \to 0$$

as $n \to \infty$. So (4.7) follows. □

We need the following two lemmas for the proof of Lemma 3.2.

LEMMA 4.3. *Let $\mathbf{\Delta}_j$ be as in (3.3) and $n_\alpha$ in (3.8). Write $\mathbf{\Delta}_j = (\Delta_{1j}, \Delta_{2j}, \ldots, \Delta_{nj})^T \in \mathbb{R}^n$. Then, for any $t > 0$,*

$$P(|\Delta_{1j}| \geq t, \ |\Delta_{2j}| \geq t) \leq e^{-t^2 n/j} + e^{-(\log n)^2/11},$$

*uniformly on $j \in (n/(\log n)^3, n_\alpha)$ as $n$ is sufficiently large.*

PROOF. Again, write $m = n_\alpha$. By (3.4), $\mathbf{\Delta}_j = \mathbf{\Gamma}_{n,j}\mathbf{\Gamma}_{n,j}^T\mathbf{y}_j$, where $\mathbf{\Gamma}_{n,j} = (\boldsymbol{\gamma}_1, \boldsymbol{\gamma}_2, \ldots, \boldsymbol{\gamma}_{j-1})$ and $\mathbf{y}_j = (y_{1j}, y_{2j}, \ldots, y_{nj})^T \in \mathbb{R}^n$. It is easy to see from the orthogonality of the $\gamma_i$'s and the independence between $\mathbf{y}_j$ and $\mathbf{\Gamma}_{n,j}$ that

(4.9) $$\mathbf{\Delta}_j \stackrel{d}{=} \mathbf{\Gamma}_{n,j}(y_{1j}, y_{2j}, \ldots, y_{j-1\,j})^T.$$

Here and later, the notation "$\stackrel{d}{=}$" means that the distributions of both sides are identical. Thus,

(4.10) $$(\Delta_{1j}, \Delta_{2j})^T \stackrel{d}{=} \left(\sum_{k=1}^{j-1}\gamma_{1k}y_{kj}, \sum_{k=1}^{j-1}\gamma_{2k}y_{kj}\right)^T.$$

Observe that $\boldsymbol{\gamma}_1, \boldsymbol{\gamma}_2, \ldots, \boldsymbol{\gamma}_{j-1}$ are functions of $\mathbf{y}_1, \mathbf{y}_2, \ldots, \mathbf{y}_{j-1}$. We know from (4.10) that $(\Delta_{1j}, \Delta_{2j})^T \sim N_2(\boldsymbol{\mu}, \boldsymbol{\Sigma})$ conditioning on $\mathbf{y}_1, \mathbf{y}_2, \ldots, \mathbf{y}_{j-1}$. Easily, $\boldsymbol{\mu} = 0$ and $\mathrm{Var}(\Delta_{pj}) \sim \sum_{k=1}^{j-1} \gamma_{pk}^2$ for $p = 1, 2$, and the correlation coefficient of $\Delta_{1j}$ and $\Delta_{2j}$ is

$$\text{(4.11)} \qquad \rho_j := \frac{\sum_{k=1}^{j-1} \gamma_{1k}\gamma_{2k}}{\sqrt{\sum_{k=1}^{j-1} \gamma_{1k}^2} \sqrt{\sum_{k=1}^{j-1} \gamma_{2k}^2}}.$$

Therefore, there exists two independent standard normals $\xi$ and $\eta$ such that the conditional distribution of $\Delta_{1j}$ and $\Delta_{2j}$ given $\mathbf{y}_1, \mathbf{y}_2, \ldots, \mathbf{y}_{j-1}$ is the same as that of $(\sum_{k=1}^{j-1} \gamma_{1k}^2)^{1/2}\xi$ and $(\sum_{k=1}^{j-1} \gamma_{2k}^2)^{1/2}(\rho_j \xi + \sqrt{1-\rho_j^2}\eta)$. It follows that

$$P(|\Delta_{1\,j+1}| \geq t,\ |\Delta_{2\,j+1}| \geq t \mid \mathbf{y}_1, \mathbf{y}_2, \ldots, \mathbf{y}_j)$$

$$\text{(4.12)} \qquad \leq P\Bigg(|\xi| \geq t\bigg(\sum_{k=1}^{j} \gamma_{1k}^2\bigg)^{-1/2},$$

$$|\eta| \geq t\bigg(\sum_{k=1}^{j} \gamma_{2k}^2\bigg)^{-1/2} - |\rho_{j+1}\xi| \ \bigg|\ \mathbf{y}_1, \mathbf{y}_2, \ldots, \mathbf{y}_j\Bigg).$$

Now, by (3.5) and (3.6), there exists a sequence of i.i.d. standard normals $\psi_1, \psi_2, \ldots, \psi_n$ such that $\mathcal{L}(\sum_{k=1}^{j} \gamma_{pk}^2) = \mathcal{L}(S_j/S_n)$ for $p = 1, 2$, where $S_j = \sum_{l=1}^{j} \psi_l^2$. By Lemma 4.2,

$$\max_{n/(\log n)^3 \leq j \leq m} P\Bigg(\bigg|\frac{S_n}{S_j} - \frac{n}{j}\bigg| \geq n^{-1/5}\Bigg)$$

$$\leq 6 \max_{n/(\log n)^3 \leq j \leq m} \Bigg\{\exp\Bigg(-\frac{j^4 n^{-2/5}}{48 n^3}\Bigg)\Bigg\} \leq e^{-\sqrt{n}},$$

as $n$ is sufficiently large. By (4.12),

$$P(|\Delta_{1\,j+1}| \geq t,\ |\Delta_{2\,j+1}| \geq t)$$

$$\text{(4.13)} \qquad \leq P(|\xi| \geq t\sqrt{(n/j) - n^{-1/5}},\ |\eta| \geq t\sqrt{(n/j) - n^{-1/5}} - |\rho_{j+1}\xi|)$$

$$+ 2e^{-\sqrt{n}}.$$

Since $P(|\xi| \geq x) \leq (1/x)\exp(-x^2/2)$ for any $x > 0$, by Lemma 4.4 below,

$$P(|\rho_{j+1}\xi| \geq (\log n)^7/n^{1/4}) \leq P\bigg(|\rho_{j+1}| \geq \frac{(\log n)^6}{\sqrt{n}}\bigg) + P(|\xi| \geq \log n)$$

$$\leq 2e^{-(\log n)^2/10}$$

for sufficiently large $n$. Thus, combining this with (4.13), we obtain from the independence of $\xi$ and $\eta$ that $P(|\Delta_{1\,j+1}| \geq t,\ |\Delta_{2\,j+1}| \geq t)$ is bounded above by

$$P(|\xi| \geq \underbrace{t\sqrt{(n/j)} - n^{-1/5}}_{A},\ |\eta| \geq \underbrace{t\sqrt{(n/j)} - n^{-1/5} - n^{-1/4}(\log n)^7}_{B})$$
$$+ 3e^{-(\log n)^2/10}$$
$$\leq 2e^{-t^2 n/j} + e^{-(\log n)^2/11},$$

uniformly on $j \in (n/(\log n)^3, m)$ as $n$ is sufficiently large, where $A$ and $B$ are essentially $t\sqrt{n/j}$ when using Lemma A.1 in the last step. $\square$

Now we measure how fast the correlation coefficient $\rho_j$ goes to zero. The idea behind the proof is that we view $\gamma_{ij}$'s in the expression of $\rho_j$ in (4.11) as independent normals with mean zero and standard deviation $n^{-1/2}$. This intuition will be carried out rigorously by using Lemma A.4.

LEMMA 4.4. *Let $\rho_j$ be as in* (4.11). *Then*

$$P(|\rho_{j+1}| \geq (\log n)^6/n^{1/4}) \leq e^{-(\log n)^2/10},$$

*uniformly on $j \in (n/(\log n)^3, n_\alpha)$ for sufficiently large $n$.*

PROOF. Write $m = n_\alpha$ for simplification. Note that $(\gamma_{11}, \gamma_{12}, \ldots, \gamma_{1n})$ has the same distribution as that of $(\gamma_{21}, \gamma_{22}, \ldots, \gamma_{2n})$ because of the Haar invariance of $\mathbf{\Gamma} = (\boldsymbol{\gamma}_1, \boldsymbol{\gamma}_2, \ldots, \boldsymbol{\gamma}_n)$. For any $a > 0$,

$$(4.14) \quad P(|\rho_{j+1}| \geq a) \leq P\left(\left|\sum_{k=1}^{j} \gamma_{1k}\gamma_{2k}\right| \geq \frac{aj}{2n}\right) + 2P\left(\left(\sum_{k=1}^{j} \gamma_{1k}^2\right)^{-1} \geq \frac{2n}{j}\right).$$

By (3.5) and (3.6), the sum appearing in the last probability in (4.14) is equal to $S_j/S_n$ in law as in Lemma 4.2. By this lemma,

$$P\left(\left|\left(\sum_{k=1}^{j} \gamma_{1k}^2\right)^{-1} - \frac{n}{j}\right| \geq \frac{n}{j}\right) = P\left(\left|\frac{S_n}{S_j} - \frac{n}{j}\right| \geq \frac{n}{j}\right)$$
$$(4.15) \qquad\qquad\qquad \leq 6\exp\left(-\frac{j^4}{48n^3}\left(\frac{n}{j}\right)^2\right)$$
$$\leq e^{-\sqrt{n}},$$

uniformly on $j \in (n/(\log n)^3, m)$ for $n$ sufficiently large. Recall (3.6) again. Choosing $m = 2,\ t = n^{-1/4}\log n, s = \log n$ and $r = (\log n)^2/\sqrt{n}$ in Theorem

A.4, by (3.6), we have $2n^2$ i.i.d. standard normals $\{y_{ij}; 1 \leq i \leq 2, \ 1 \leq j \leq n\}$ such that

$$P\left(\varepsilon_n(2) \geq \frac{(\log n)^2}{n^{1/4}}\right)$$

$$\leq 4n^2 \exp\left(-\frac{(\log n)^4}{16}\right)$$

(4.16)

$$+ 3n^2 e^{-(\log n)^2/2} + 3n^{5/4}\left(1 + \frac{(\log n)^2}{3\sqrt{n}(\sqrt{n}+2)}\right)^{-n/2}$$

$$\leq e^{-(\log n)^2/9}$$

for $n$ large enough, where $\varepsilon_n(2) = \max_{1 \leq i \leq 2, 1 \leq j \leq n} |\sqrt{n}\gamma_{ij} - y_{ij}|$. Notice that

(4.17) $\quad n\left|\sum_{k=1}^{j} \gamma_{1k}\gamma_{2k}\right| \leq \left|\sum_{k=1}^{j} y_{1k}y_{2k}\right| + \frac{(\log n)^2}{n^{1/4}} \sum_{i=1}^{2} \sum_{k=1}^{j} |y_{ik}| + \frac{2j(\log n)^4}{\sqrt{n}}$

on $\{\varepsilon_n(2) \leq (\log n)^2/n^{1/4}\}$. Note that $E\exp(|y_{11}y_{21}|/8) < \infty$ and $E|y_{11}| \leq 1$. By Lemma A.3, there exists a universal constant $C > 0$ such that

(4.18)
$$P\left(\sum_{i=1}^{2}\sum_{k=1}^{j} |y_{ik}| \geq 3j\right) \leq e^{-Cj} \quad \text{and}$$

$$P\left(\left|\sum_{k=1}^{j} y_{1k}y_{2k}\right| \geq \sqrt{j}\log j\right) \leq e^{-(\log n)^2/3},$$

uniformly on $j \in (n/(\log n)^3, m)$, where the first one comes from (i) of Lemma A.3 and the second is obtained by (ii) of Lemma A.3 in the same way as in (4.8). If neither of the events in the above two probabilities occurs and $\varepsilon_n(2) \leq (\log n)^2/n^{1/4}$, then from (4.17)

$$n\left|\sum_{k=1}^{j} \gamma_{1k}\gamma_{2k}\right| \leq \sqrt{j}\log j + \frac{3j(\log n)^2}{n^{1/4}} + \frac{2j(\log n)^4}{\sqrt{n}} < 5n^{3/4}(\log n)^2,$$

uniformly on $j \in (n/(\log n)^3, m)$ for sufficiently large $n$. Thus, from (4.16) and (4.18),

$$P\left(\left|\sum_{k=1}^{j} \gamma_{1k}\gamma_{2k}\right| \geq \frac{5(\log n)^2}{n^{1/4}}\right) \leq 2e^{-(\log n)^2/9},$$

as $n$ is sufficiently large. Choose $a = (\log n)^6/n^{1/4}$ in (4.14). Then, $aj/(2n) \geq 5(\log n)^2/n^{1/4}$ for all $j \in (n(\log n)^{-3}, m)$, as $n$ is sufficiently large. It follows from the above that

(4.19) $\quad P\left(\left|\sum_{k=1}^{j} \gamma_{1k}\gamma_{2k}\right| \geq \frac{aj}{2n}\right) \leq 2e^{-(\log n)^2/9},$

uniformly on $j \in (n/(\log n)^3, m)$ as $n$ is sufficiently large. It is easy to see that the last probability in (4.14) is bounded by the first probability in (4.15). Combining (4.14), (4.15) and (4.19) together, we obtain that

$$P(|\rho_{j+1}| \geq (\log n)^6/n^{1/4}) \leq 2e^{-(\log n)^2/9} + 2e^{-\sqrt{n}} \leq e^{-(\log n)^2/10},$$

as $n$ is sufficiently large. $\square$

PROOF OF LEMMA 3.2. Write $m = n_\alpha$. Rewrite $\boldsymbol{\Delta}_{k+1} = (\Delta_{1,k+1}, \Delta_{2,k+1}, \ldots, \Delta_{n,k+1})^T \in \mathbb{R}^n$. By (4.9) and (4.10), $\mathcal{L}(\Delta_{i,k+1}) = \mathcal{L}(\sum_{l=1}^{k} \gamma_{il} y_{lk+1})$, so conditioning on $\mathbf{y}_1, \mathbf{y}_2, \ldots, \mathbf{y}_k$,

$$(4.20) \qquad \Delta_{i,k+1} \sim N\left(0, \sum_{l=1}^{k} \gamma_{il}^2\right).$$

Let $E_0 = \{\max_{2 \leq j \leq k} \|\boldsymbol{\Delta}_j\| \leq t\}$ and $E_i = \{|\Delta_{i,k+1}| \leq t\}$. Although each $E_i$ depends on $n$ and $k$, we would rather use the notation $E_i$ for simplification. This will not cause confusion in the context. Evidently,

$$(4.21) \qquad \left\{\max_{2 \leq j \leq k+1} \|\boldsymbol{\Delta}_j\| \leq t\right\} = \bigcap_{i=0}^{n} E_i.$$

To apply Lemma 4.1, we now calculate $P(E_0 \backslash E_i)$. Define

$$\delta_n = \max_{(i,l) \in \Omega_n} \left|\left(\sum_{j=1}^{l} \gamma_{ij}^2\right)^{-1/2} - \sqrt{\frac{n}{l}}\right|,$$

where

$$\Omega_n = \{(i,l);\ 1 \leq i \leq n,\ n/(\log n)^3 \leq l \leq m\}.$$

Recall (4.20). Let $S_j$ be as in Lemma 4.2, then by the lemma and the fact that $|\sqrt{a} - \sqrt{b}| \leq |a - b|$ if $a \geq 1$,

$$(4.22)\quad P\left(\delta_n \geq \frac{(\log n)^8}{\sqrt{n}}\right) \leq n^2 \max P\left(\left|\frac{S_n}{S_l} - \frac{n}{l}\right| \geq \frac{(\log n)^8}{\sqrt{n}}\right) \leq e^{-(\log n)^2}$$

for sufficiently large $n$, where the max above is taken over all $l$ such that $n/(\log n)^3 \leq l \leq m$. By (4.20), for some standard normal $\xi$, we have $\Delta_{i,k+1} \sim (\sum_{j=1}^{k} \gamma_{ij}^2)^{1/2}\xi$ conditioning on $\mathbf{y}_1, \mathbf{y}_2, \ldots, \mathbf{y}_k$. Thus, $P(E_i^c \mid \mathbf{y}_1, \mathbf{y}_2, \ldots, \mathbf{y}_k) = P(|\xi| > (\sum_{j=1}^{k} \gamma_{ij}^2)^{-1/2} t \mid \mathbf{y}_1, \mathbf{y}_2, \ldots, \mathbf{y}_k)$. It follows that on $\{\delta_n \leq (\log n)^8/\sqrt{n}\}$,

$$(4.23)\qquad \begin{aligned} f_n^+(k) &= P\left(|\xi| > t\left(\sqrt{\frac{n}{k}} + \frac{(\log n)^8}{\sqrt{n}}\right)\right) \\ &\leq P(E_i^c \mid \mathbf{y}_1, \mathbf{y}_2, \ldots, \mathbf{y}_k) \\ &\leq P\left(|\xi| > t\left(\sqrt{\frac{n}{k}} - \frac{(\log n)^8}{\sqrt{n}}\right)\right) = f_n^-(k), \end{aligned}$$

uniformly on $(i,k) \in \Omega_n$. The key observation for this proof is that the above conditional probability is bounded above and below by unconditional probabilities. Obviously, $E_0$ is a set in the $\sigma$-algebra generated by $\mathbf{y}_1, \mathbf{y}_2, \ldots, \mathbf{y}_k$. By (4.22) and (4.23),

$$P(E_0 \backslash E_i) = E\{I_{E_0}(P(E_i^c \mid \mathbf{y}_1, \mathbf{y}_2, \ldots, \mathbf{y}_k))\} \leq P(E_0) f_n^-(k) + e^{-(\log n)^2}$$

for all $(i,k) \in \Omega_n$ when $n$ is sufficiently large. Similarly, use the first step above to obtain

$$P(E_0 \backslash E_i) \geq P(E_0 \cap F_n) \cdot f_n^+(k) \geq P(E_0) \cdot f_n^+(k) - e^{-(\log n)^2},$$

for all $(i,k) \in \Omega_n$, where $F_n := \{\delta_n \leq (\log n)^8/\sqrt{n}\}$. Therefore,

$$(4.24) \quad nP(E_0) \cdot f_n^+(k) - ne^{-(\log n)^2} \leq \sum_{i=1}^n P(E_0 \backslash E_i)$$
$$\leq nP(E_0) \cdot f_n^-(k) + ne^{-(\log n)^2},$$

uniformly on $n/(\log n)^3 \leq k \leq m$ as $n$ is sufficiently large.

Finally, note that $e^{-t^2 n/j}$ is increasing in $j$. By Lemma 4.3, $P(E_1^c E_2^c) \leq n^{-t^2/\alpha}(\log n)^C$ for some constant $C > 0$ as $n$ is sufficiently large. Also, the $n$ random variables in $(\Delta_{1,k+1}, \Delta_{2,k+1}, \ldots, \Delta_{n,k+1})$ are exchangeable by the Haar-invariance. Hence,

$$(4.25) \quad \sum_{1 \leq i < j \leq n} P(E_i^c E_j^c) \leq \frac{n^2}{2} P(E_1^c E_2^c) \leq \frac{(\log n)^C}{n^{t^2/\alpha - 2}},$$

as $n$ is sufficiently large. By (4.24), the quantity $P(E_0) - \sum_{i=1}^n P(E_0 \backslash E_i)$ is bounded above and below respectively by

$$(1 - nf_n^+(k))P(E_0) + ne^{-(\log n)^2} \text{ and } (1 - nf_n^-(k))P(E_0) - ne^{-(\log n)^2}.$$

This together with (4.25) yields the desired conclusion via Lemma 4.1. $\square$

## APPENDIX

The following is a standard result. It can be found in, for example, Lemma 3 on page 49 from [6].

LEMMA A.1. *Suppose $X \sim N(0,1)$. Then*

$$\frac{1}{\sqrt{2\pi}} \cdot \frac{x}{1+x^2} e^{-x^2/2} \leq P(X > x) \leq \frac{1}{\sqrt{2\pi}} \cdot \frac{1}{x} e^{-x^2/2}$$

*for all $x > 0$.*

The following lemma is part (ii) on page 186 from [29].

LEMMA A.2. *Suppose $\mathbf{y}$ is an $\mathbb{R}^n$-valued random vector with multinormal distribution with mean $\mathbf{0}$ and covariance matrix $\mathbf{\Sigma}$ of rank $r$. If $\mathbf{\Sigma}^2 = \mathbf{\Sigma}$, then there exists a sequence of independent standard normals $\{\xi_j;\ j = 1, 2, \ldots, n\}$ such that $\|\mathbf{y}\|^2$ has the same distribution as that of $\sum_{j=1}^r \xi_j^2$, that is, $\|\mathbf{y}\|^2 \sim \chi^2(r)$.*

For $A \subset \mathbb{R}$, the interior and the closure of $A$ in $\mathbb{R}$ are denoted by $A^\circ$ and $\bar{A}$, respectively. The following are Chernoff's bound and a moderate deviation result. They can be found from, for example, (c) of Remarks on page 27 from [9] and Theorem 3.7.1 on page 109 from [9].

LEMMA A.3. *Let $\{X, X_i, i = 1, 2, \ldots\}$ be a sequence of i.i.d. random variables. Let $S_n = \sum_{i=1}^n X_i, n \geq 1$. Then:*
(i) *For any $A \subset \mathbb{R}$ and $n \geq 1$,*

$$P(S_n/n \in A) \leq 2e^{-nI(A)},$$

*where $I(x) = \sup_{t \in \mathbb{R}}\{tx - \log E(e^{tX})\}$ and $I(A) = \inf_{x \in A} I(x)$.*
(ii) *Assume further that $EX = 0$, $\operatorname{var}(X) = \sigma^2 > 0$ and $Ee^{t_0 X} < \infty$ for some $t_0 > 0$. Let $\{a_n; n = 1, 2, \ldots\}$ be a sequence of positive numbers such that $a_n \to 0$ and $na_n \to \infty$ as $n \to \infty$. Then*

$$\lim_{n \to \infty} a_n \log P\left(\sqrt{\frac{a_n}{n}} S_n \in A\right) = -\inf_{x \in A}\left\{\frac{x^2}{2\sigma^2}\right\}$$

*for any subset $A \subset \mathbb{R}$ such that $\inf\{|x|; x \in A^\circ\} = \inf\{|x|; x \in \bar{A}\}$.*

The following lemma is Theorem 5 from [19].

LEMMA A.4. *For each $n \geq 2$, there exists matrices $\mathbf{\Gamma}_n = (\gamma_{ij})_{1 \leq i, j \leq n}$ and $\mathbf{Y}_n = (y_{ij})_{1 \leq i, j \leq n}$ whose $2n^2$ elements are random variables defined on the same probability space such that:*
(i) *the law of $\mathbf{\Gamma}_n$ is the normalized Haar measure on the orthogonal group $O_n$;*
(ii) *$\{y_{ij}; 1 \leq i, j \leq n\}$ are i.i.d. random variables with the standard normal distribution;*
(iii) *set $\varepsilon_n(m) = \max_{1 \leq i \leq n, 1 \leq j \leq m} |\sqrt{n}\gamma_{ij} - y_{ij}|$ for $m = 1, 2, \ldots, n$. Then*

$$P(\varepsilon_n(m) \geq rs + 2t) \leq 4me^{-nr^2/16} + 3mn\left(\frac{1}{s}e^{-s^2/2} + \frac{1}{t}\left(1 + \frac{t^2}{3(m + \sqrt{n})}\right)^{-n/2}\right)$$

*for any $r \in (0, 1/4)$, $s > 0$, $t > 0$, and $m \leq (r/2)n$.*

**Acknowledgments.** The author very much thanks Amir Dembo, Persi Diaconis, Morris Eaton, William Sudderth and Ofer Zeitouni for very useful discussions, email and correspondence.

# REFERENCES

[1] ANDERSON, T. W. (1984). *An Introduction to Multivariate Statistical Analysis*, 2nd ed. Wiley, New York. MR0771294

[2] APOSTOL, T. M. (1974). *Mathematical Analysis*, 2nd ed. Addison–Wesley, Reading, MA. MR0344384

[3] BAI, Z. D. (1999). Methodologies in spectral analysis of large-dimensional random matrices, a review. *Statist. Sinica* **9** 611–677. MR1711663

[4] BILLINGSLEY, P. (1979). *Probability and Measure*. Wiley, New York. MR0534323

[5] BOREL, E. (1906). *Introduction géometrique á quelques théories physiques*. Gauthier–Villars, Paris. JFM 45.0808.10

[6] CHOW, Y. S. and TEICHER, H. (1988). *Probability Theory, Independence, Interchangeability, Martingales*, 2nd ed. Springer, New York. MR0953964

[7] COLLINS, B. (2003). Intégrales matricielles et probabilitiés non-commutatives. Thèse de Doctorat, Univ. Paris 6.

[8] D'ARISTOTLE, A., DIACONIS, P. and NEWMAN, C. M. (2002). Brownian motion and the classical groups. *Probability, Statistics and Their Applications: Papers in Honor of Rabi Bhattacharya* 97–116. IMS Lecture Notes Monogr. Ser. **41**. IMS, Beachwood, OH. MR1999417

[9] DEMBO, A. and ZEITOUNI, O. (1998). *Large Deviations Techniques and Applications*, 2nd ed. Springer, New York. MR1619036

[10] DIACONIS, P. (2003). Patterns in eigenvalues: The 70th Josiah Willard Gibbs lecture. *Bull. Amer. Math. Soc.* (*N.S.*) **40** 155–178 (electronic). MR1962294

[11] DIACONIS, P. and FREEDMAN, D. (1987). A dozen de Finetti-style results in search of a theory. *Ann. Inst. H. Poincaré Probab. Statist.* **23** 397–423. MR0898502

[12] DIACONIS, P. and SHAHSHAHNI, M. (1994). On the eigenvalues of random matrices. Studies in applied probability. *J. Appl. Probab.* **31A** 49–62. MR1274717

[13] DIACONIS, P. and EVANS, S. N. (2001). Linear functionals of eigenvalues of random matrices. *Trans. Amer. Math. Soc.* **353** 2615–2633 (electronic). MR1828463

[14] DIACONIS, P. W., EATON, M. L. and LAURITZEN, S. L. (1992). Finite de Finetti theorems in linear models and multivariate analysis. *Scand. J. Statist.* **19** 289–315. MR1211786

[15] EATON, M. L. (1989). *Group Invariance Applications in Statistics*. IMS, Hayward, CA. MR1089423

[16] GALLARDO, L. (1983). Au sujet du contenu probabiliste d'un lemma d'Henri Poincaré. *Ann. Univ. Clemont* **69** 192–197. MR0645576

[17] GEMAN, S. (1980). A limit theorem for the norm of random matrices. *Ann. Probab.* **8** 252–261. MR0566592

[18] HORN, R. and JOHNSON, C. (1990). *Matrix Analysis*. Cambridge Univ. Press. MR1084815

[19] JIANG, T. (2005). Maxima of entries of Haar distributed matrices. *Probab. Theory Related Fields* **131** 121–144. MR2105046

[20] JOHANSSON, K. (1997). On random matrices from the compact classical groups. *Ann. of Math.* (*2*) **145** 519–545. MR1454702

[21] JONSSON, D. (1982). Some limit theorems for the eigenvalues of a sample covariance matrix. *J. Multivariate Anal.* **12** 1–38. MR0650926

[22] LANG, S. (1987). *Calculus of Several Variables*. Springer, New York. MR0886677

[23] LETAC, G. (1981). Isotropy and sphericity: Some characterisations of the normal distribution. *Ann. Statist.* **9** 408–417. MR0606624

[24] MAXWELL, J. C. (1875). *Theory of Heat*, 4th ed. Longmans, London.

[25] MAXWELL, J. C. (1878). On Boltzmann's theorem on the average distribution of energy in a system of material points. *Cambridge Phil. Soc. Trans.* **12** 547. JFM 11.0776.01
[26] MCKEAN, H. P. (1973). Geometry of differential space. *Ann. Probab.* **1** 197–206. MR0353471
[27] POINCARÉ, H. (1912). *Calcul des probabilitiés*. Gauthier–Villars, Paris. JFM 43.0308.04
[28] RAINS, E. M. (1997). High powers of random elements of compact Lie groups. *Probab. Theory Related Fields* **107** 219–241. MR1431220
[29] RAO, C. R. (1973). *Linear Statistical Inference and Its Applications*. Wiley, New York. MR0346957
[30] STAM, A. J. (1982). Limit theorems for uniform distributions on spheres in high-dimensional Euclidean spaces. *J. Appl. Probab.* **19** 221–228. MR0644435
[31] YIN, Y. Q., BAI, Z. D. and KRISHNAIAH, P. R. (1988). On the limit of the largest eigenvalue of the large-dimensional sample covariance matrix. *Probab. Theory Related Fields* **78** 509–521. MR0950344
[32] YOR, M. (1985). *Inégalitiés de martingales continus arrêtès à un temps quelconques* I. *Lecture Notes in Math.* **1118**. Springer, Berlin. Zbl 0563.60045

SCHOOL OF STATISTICS
UNIVERSITY OF MINNESOTA
313 FORD HALL
224 CHURCH STREET S.E.
MINNEAPOLIS, MINNESOTA 55455
USA
E-MAIL: tjiang@stat.umn.edu